\documentclass[11pt, twoside,leqno]{amsart}
\usepackage{geometry}
\allowdisplaybreaks[1]  
\usepackage[utf8]{inputenc}
\usepackage{amsfonts, amssymb, amscd, amsmath, amsthm}
\usepackage{setspace}
\usepackage[utf8]{inputenc}
\usepackage{graphicx}
 \usepackage{float}
\usepackage{fullpage}
\usepackage{mathtools}
\newtheorem{thm}{Theorem}[section]
\newtheorem{lemma}[thm]{Lemma}
\newtheorem{prop}[thm]{Proposition}
\newtheorem{cor}[thm]{Corollary}
\newtheorem{defn}[thm]{Definition}

\newtheorem{note}[thm]{Note}

\newcommand{\B}{\mathcal{B}}
\newcommand{\N}{\mathbb{N}}

\newcommand{\A}{\mathcal{A}}

 \usepackage[bookmarksnumbered, colorlinks, plainpages]{hyperref}
  \hypersetup{colorlinks=true,linkcolor=blue, anchorcolor=green, citecolor=cyan, urlcolor=blue, filecolor=magenta, pdftoolbar=true}

\DeclareMathOperator{\Dom}{Dom}
\DeclareMathOperator{\Span}{Span}

\usepackage[charter]{mathdesign}
\usepackage{comment}
\title{A Radon-Nikod\'ym theorem for local completely positive invariant multilinear maps}
\author{Anindya Ghatak and Santhosh Kumar Pamula}
\address{Statistics and Mathematics Unit, Indian Statistical Institute, R. V. College Post, Bangalore 560059, India}
\email{anindya.ghatak123@gmail.com}

\address{Department of Mathematics, Birla Institute of Technology and Science (BITS Pilani), Pilani
campus, Rajasthan 333031, India.}
\email{santhoshk.pamula@pilani.bits-pilani.ac.in}
\subjclass[2010]{46L07, 46L05, 46G25, 47L40}
\keywords{Locally $C^{\ast}$-algebra, Multilinear map, Local completely positive map, Stinespring's theorem, Radon-Nikod\'{y}m theorem.}
\date{\today}
\begin{document}
\maketitle 
\begin{abstract}
	In this article, we introduce local completely positive $k$-linear maps between locally $C^{\ast}$-algebras and obtain  Stinespring type representation by adopting the notion of ``invariance'' defined by J. Heo for $k$-linear maps between $C^{\ast}$-algebras. Also, we supply the minimality condition to make certain that minimal representation is unique up to unitary equivalence.   As a consequence, we prove Radon-Nikod\'{y}m theorem for unbounded operator-valued local completely positive invariant $k$-linear maps. The obtained Radon-Nikod\'{y}m derivative is a positive contraction on some Hilbert space with several reducing subspaces.  
\end{abstract}
 

\section{Introduction}

In 1987, Christensen and Sinclair introduced multilinear completely positive (in short CP) maps  \cite{Christensen 87} between $C^*$-algebras and used this notion to develop theory on the cohomology group of the von Neuman algebras. These authors gave a representation for completely bounded symmetric multilinear maps in terms of $\ast$-representations 
of the algebra and constructed suitable bridging operators between the Hilbert spaces by using the famous result of Wittstock's matricial Hahn Banach theorem \cite{Wittstock 83}. In the article \cite{SS98}, Sinclair and Smith proved that a von Neuman algebra is injective if and only if every bilinear completely bounded (CB) maps can be factorized in terms of bilinear CP-maps. Paulsen and Smith \cite{Paulsen 87} extended this result to operator spaces by using Wittstock's matricial Hahn-Banach theorem. Further, with the help of the representation of completely bounded symmetric multilinear maps, Effros et. al. showed in  \cite{Effros 87} that the completely bounded cohomology group from a $C^{\ast}$-algebra to an injective von Neumann algebra is zero.
  

    To study the structure of multilinear CP-maps, Christensen established Stinespring type theorem for symmetric $k$-linear CP-maps (see \cite[Section 3.]{Christensen 87} for details).
     Motivated by this work,  Heo \cite{Heo 00} came up with an additional property (called as {\it \bf  invariant}) on  $k$-linear CP-maps as an analogue of invariant sesquilinear forms \cite{AI 90}. With this additional invariance property, Heo obtained a dilation for invariant $k$-linear CP-map to a multilinear representation which can be factorized as a product of $m$-number of commuting representations, where $m= [\frac{k+1}{2}].$ Further, the author developed a Radon-Nikod\'{y}m theorem for completely bounded and completely positive invariant $k$-linear maps. Several authors studied invariant multilinear maps, for example, see \cite{GK09, Gudder79, HJ 2019, Heo 2016, Heo10, KRP96} and references therein. In this article, our main concern is to study the structure of $k$-linear local CP-maps between locally $C^*$-algebras. These results generalize the work of Heo \cite{Heo}.

  We organize this manuscript into five sections. In the second section, we fix notations, recall some definitions and known results from the literature which are needed for the following sections. In the third section, firstly we construct a new quantized domain (contains the original quantized domain) associated to a given local completely positive invariant $k$-linear map defined from a locally $C^{\ast}$-algebra $\mathcal{A}$ to a locally $C^{\ast}$-algebra of non-commutative continuous functions on a quantized domain. Next, we construct a commuting family of local contractive $\ast$-homomorphisms from $\mathcal{A}$ to a locally $C^{\ast}$-algebra of non-commutative continuous functions on the new quantized domain. Finally, we show that the product of these local contractive $\ast$-homomorphisms is a desired dilation of the given invariant local completely positive $k$-linear map (see Theorem \ref{Theorem: Multilinear Stinespring}). This result extends the Anar Dosiev's Stinespring theorem \cite[Theorem 5.1]{Dosiev} to $k$-linear maps and a locally convex analogue of Heo's result \cite{Heo}.
 
 In the fourth section, a suitable notion of minimality for our Stinespring type representation has been identified to ensure the uniqueness of minimal representation up to unitary equivalence (see Theorem \ref{Theorem:unitary equivalence of minimal Stinespring representation}). In the final section, we prove Radon-Nikod\'{y}m theorem for multilinear local CP maps that characterizes the set of all maps dominated by a fixed local completely positive invariant $k$-linear map. Here, the obtained Radon-Nikod\'{y}m derivative is a positive contraction on some Hilbert space with an upward filtered family of reducing sub spaces (see Theorem \eqref{Theorem: Radon-Nikodym}).
\section{Preliminaries}
A locally $C^*$-algebra $\mathcal{A}$ is a $*$-algebra which is complete in a locally convex topology determined by an upward filtered family $\{p_{\alpha}: \alpha\in \Lambda\}$ of $C^*$-seminorms. In this article, 
we denote $\mathcal{A}$ and $\mathcal{B}$ as the unital locally $C^*$-algebras with the upward filtered family of $C^*$-seminorms  $\{p_{\alpha}: \alpha\in \Lambda\}$ and $\{q_{\ell}: \ell\in \Omega\}$ on $\mathcal{A}$ and
$\mathcal{B}$ respectively. Locally $C^*$-algebras are the generalization of $C^*$-algebras, indeed it can be seen as the  projective limit of the inverse system of $C^{*}$-algebras.
For instance, let $I_{\alpha}=\{ a\in \mathcal{A}: p_{\alpha}(a)=0\}$. Then $I_{\alpha}$ is a closed two-sided ideal in $\mathcal{A}$ and hence  $\mathcal{A}_{\alpha}:= \mathcal{A}/I_{\alpha}$  is a $C^*$-algebra with the $C^{\ast}$-norm $\|\cdot \|_{\alpha}$ given by
\begin{equation*}
\|a_{\alpha}\|_{\alpha} := p_{\alpha}(a), \; \text{for all}\; a_{\alpha} := a+I_{\alpha}\in \mathcal{A}_{\alpha},\; \alpha \in \Lambda.
\end{equation*} 
For every $\alpha, \beta \in \Lambda$ there exist a canonical projection $\pi_{\alpha} \colon \mathcal{A} \to \mathcal{A}_{\alpha}$ and a canonical $\ast$-epimorphism $\pi_{\alpha, \beta} \colon \mathcal{A}_{\beta} \to \mathcal{A}_{\beta}$ whenever $\alpha \leq \beta$ satisfying the relation $\pi_{\alpha} = \pi_{\alpha, \beta} \circ \pi_{\beta}$. The projective limit of an inverse system $\big( \{\mathcal{A}_{\alpha}\}_{\alpha \in \Lambda}; \{\pi_{\alpha, \beta}\}_{\alpha \leq \beta}\big)$ of $C^{\ast}$-algebras is defined as 
\begin{equation*}
 \lim\limits_{\underset{\alpha}{\longleftarrow}} \mathcal{A}_{\alpha}:=\Big\{\{a_{\alpha}\}_{\alpha\in \Lambda}  \in \prod_{\alpha\in \Lambda} \mathcal{A}_{\alpha}: \pi_{\alpha\beta}(a_{\beta})= 
 a_{\alpha} \text{~~whenever~~} \alpha \leq \beta,\; \alpha, \beta\in \Lambda \Big\}.
 \end{equation*}
 Then $\mathcal{A}$ is identified with $\lim\limits_{\underset{\alpha}{\longleftarrow}} \mathcal{A}_{\alpha}$ via the following map
 \begin{align}\label{Eq: Identification by projective limits}
     a \to \{\pi_{\alpha}(a)\}_{\alpha \in \Lambda} \; \text{for all}\;  a\in \mathcal{A}.
\end{align}
 From the identification given in Equation \eqref{Eq: Identification by projective limits}, it is natural to see the existence of local self-adjoints and local positive elements in local $C^*$-algebras.  For further details, we refer the reader to \cite{Gheo 17}.
\begin{defn}\label{Defn: local self-adjoint, positive}\cite{Dosiev} An element $a\in \mathcal{A}$ is called
\begin{enumerate}
    \item  local self-adjoint if $a^*=a+x$, where $x\in \mathcal{A}$ such that $p_{\alpha}(x)=0$ for some $\alpha \in \Lambda;$
    \item local positive if $a= b^*b+y,$ where $b, y\in \mathcal{A}$ such that  $p_{\beta}(y)=0$ for some $\beta \in \Lambda.$
\end{enumerate}
\end{defn}
In other words, $a$ is local self-adjoint (local positive) if and only if $\pi_{\alpha}(a)$ is self-adjoint (positive) in $\mathcal{A}_{\alpha},$ for some $\alpha \in \Lambda.$ For brevity, we write $a^*=_{\alpha}a $ if $\pi_{\alpha}(a)$ is self-adjoint in $\mathcal{A}_{\alpha}$ and $a \geq_{\alpha} 0$ if $\pi_{\alpha}(a)$ positive in $\mathcal{A}_{\alpha}.$ In 2008, Anar Dosi defined the notion of local completely positive maps between local operator systems, which is a locally convex analogues of operator systems. Initially, such maps between local operator spaces have been studied by Effros (see \cite{Eff-Web-96} for details). 

Now we turn our discussion to multilinear maps between locally $C^{\ast}$-algebras. Let $\mathcal{A}^k= \mathcal{A}\times
\hdots\times\mathcal{A}$ denote the Cartesian product of $k$-copies of $\A,$ where $k \in \mathbb{N}.$ Then
 $$\mathcal{A}^k:= \big\{(a_{1}, \hdots, a_{k}): a_{i}\in \mathcal{A}, 1\leq i\leq k\big\}$$ is a locally convex space with associated seminorm $\Vert \cdot \Vert_{\alpha}$ for each $\alpha\in \Lambda$ given by 
\begin{equation*}
    \Vert(a_{1}, \hdots,  a_{k})\Vert_{\alpha}= \max\{ p_{\alpha}(a_{i}): 1\leq i\leq k\},\; \text{for all}\; (a_{1}, \hdots,  a_{k})\; \in \mathcal{A}^{k}.
\end{equation*}
\begin{defn}\label{Definition: alpha-symmetric}($\alpha$-symmetric)
An element $(a_{1}, a_{2}, \hdots, a_{k}) \in \mathcal{A}^{k}$ is said to be $\alpha$-symmetric (for $\alpha \in \Lambda$) if $a_{p}^{\ast} =_{\alpha}a_{k+1-p} $ for all $p=1,2,\cdots,m$ where $m = \big[\frac{k+1}{2} \big].$\\
\end{defn}
Now we define the notion of symmetric multilinear maps in local set up.
A map $\varphi: \mathcal{A}^{k}\to \mathcal{B}$
is called \emph{$k$-linear map} if 
\begin{align}
    \varphi(a_{1}, \hdots, a_{i}+ \lambda b_{i}, \hdots, a_{k})= \varphi(a_{1}, \hdots, a_{i}, \hdots, a_{k})+ \lambda\varphi(a_{1}, \hdots, b_{i}, \hdots, a_{k}),
\end{align}
for each $a_{i}, b_{i} \in \mathcal{A}, \; \lambda \in \mathbb{C}$ and $1\leq i \leq k.$
 The adjoint $\varphi^{\ast}$ of a $k$-linear map $\varphi$ is  defined as,
\begin{equation}
     \varphi^*(a_{1}, a_{2}, a_{3}, \hdots, a_{k})= \varphi(a_{k}^*, \hdots,a_{3}^{\ast}, a_{2}^*, a_{1}^*)^*, \text{~for all~} a_{1}, \hdots , a_{k}\in \mathcal{A}.
 \end{equation}
A $k$-linear map is called {\it symmetric}, if $\varphi = \varphi^{\ast}$. In other words, $\varphi$ is symmetric if and only if 
$$
\varphi(a_{1}, a_{2}, a_{3}, \hdots, a_{k})^*=\varphi(a_{k}^*, \hdots,a_{3}^{\ast}, a_{2}^*, a_{1}^*), \;  \text{for all}\; a_{1}, \hdots , a_{k}\in \mathcal{A}.
$$
Since $\mathcal{A}$ is a locally $C^{\ast}$-algebra, the set of all $n\times n$ matrices over $\mathcal{A}$ denoted by $M_{n}(\mathcal{A})$ is a local $C^*$-algebra with an upward filtered family $\{p_{\alpha}^{n}:\; \alpha \in \Lambda\}$ of $C^{\ast}$-seminorms.
 \begin{defn}\cite{Heo}\label{Definition: amplification}
     For every $n \in \mathbb{N}$, the amplification map $\varphi_{n} \colon M_{n}(\mathcal{A})^k \to M_{n}(\mathcal{B})$ is defined by 
 \begin{align*}
        \varphi_{n}(A_{1}, A_{2}, &\hdots, A_{k})
        = \Bigg[
         \sum\limits_{r_{1}, r_{2}, \hdots, r_{k-1} = 1}^{n} \varphi\big(a_{1_{ir_{1}}}, a_{2_{r_{1}r_{2}}}, \hdots, a_{k_{r_{(k-1)}j}}\big)\Bigg]_{i,j=1}^{n},
    \end{align*}
    where $A_{t} = \big[a_{t_{ij}}\big]_{i,j=1}^{n} \in M_{n}(\mathcal{A})$, for $t = 1, \hdots, k$.
 \end{defn}
 Now, we introduce  local CP-maps and  local CB-maps in this context (cf. \cite{Paulsen} for CP-maps).  Throughout this article, we write $m = \big[\frac{k+1}{2}\big]$, where $\big[ \cdot \big ]$ indicates the greatest integer function.
\begin{defn}\label{Def: LCC} Let $\varphi: \mathcal{A}^{k}\to  \mathcal{B}$  be a $k$-linear map. Then $\varphi$ is said to be 
\begin{enumerate}
\item \emph{local contractive,} if for each $\ell\in \Omega,$ there exists an $\alpha_{\ell}\in \Lambda$ such that 
\begin{equation*}
q_{\ell}(\varphi(a_{1}, a_{2}, \hdots, a_{k})) \leq p_{\alpha_{\ell}} ((a_{1}, a_{2}, \hdots ,a_{k})), \; \text{for all }\; (a_{1}, a_{2}, \hdots, a_{k})\in \mathcal{A}^{k}.
\end{equation*}
\item \emph{local completely contractive,} if the map $\varphi_{n}$ is local contractive for each $n \in \mathbb{N}$ with an additional condition that given $\ell\in \Omega,$ the existence of $\alpha_{\ell}\in \Lambda$ is independent of $n$. \end{enumerate}
\end{defn}
\begin{defn}\label{Def: LCP}    
    Let $\varphi: \mathcal{A}^{k}\to  \mathcal{B}$  be a $k$-linear map and let $m = \big[\frac{k+1}{2}\big]$. Then $\varphi$ is said to be
     \begin{enumerate}
\item \emph{local positive,} if for every $\ell \in \Omega$, there exists an $\alpha_{\ell} \in \Lambda$  such that
\begin{align*}
     & \varphi(a_{1}, a_{2}, \hdots, a_{k})\geq_{l} 0, \; \text{whenever}\; a_{p}^{\ast} =_{\alpha_{\ell}} a_{k+1-p}\; \text{for all}\; p\in \{1,2,3,\hdots, m\} \text{ and }\\
    &\varphi(a_{1}, a_{2}, \hdots, a_{k}) =_{l} 0 \text{~whenever~} a_{i}=_{\alpha_{\ell}}0 \text{ for some }1\leq i\leq k, \text{and additionally}\\
   & a_{m} \geq_{\alpha_{\ell}} 0 \text{ if } k \text{ is odd.}
    \end{align*}
 \item \emph{local completely positive,} if  the map  $\varphi_{n}$ is local positive for each $n \in \mathbb{N}$ with an additional condition that given $\ell\in \Omega,$ the existence of $\alpha_{\ell}\in \Lambda$ is independent of $n$.
 \end{enumerate}
 \end{defn}
Note that if $k=1,$ then  Definitions \ref{Def: LCC},\;\ref{Def: LCP}  coincide with the notion of local completely contractive and local completely positive linear maps defined  between local $C^*$-algebras \cite{Dosiev}. Furthermore, if $\A, \B$ are $C^*$-algebras, then Definitions \ref{Def: LCC},\;\ref{Def: LCP}  coincide with the notion of  completely contractive and completely positive $k$-linear maps studied in \cite{Heo 00, Heo}. Now we introduce the notion of invariance from \cite{Heo} to  $k$-linear maps on locally $C^*$-algebras.
 \begin{defn}
A $k$-linear map $\varphi: \mathcal{A}^{k}\to \mathcal{B}$ is said to be \emph{invariant} if: 
\begin{enumerate}
\item for $k=2m-1$ (odd), we have
\begin{align*}
    \varphi(a_{1}c_{1},\hdots, a_{m-1} c_{m-1},& a_{m}, a_{m+1}, \hdots, a_{2m-1})\\
    &=\varphi(a_{1}, \hdots, a_{m-1}, a_{m}, c_{m-1}a_{m+1}, \hdots, c_{1}a_{2m-1}),
\end{align*} 
for all $a_{1}, a_{2}, \hdots, a_{2m-1}, c_{1}, \hdots, c_{m-1} \in \mathcal{A}$.
\item  for $k=2m$ (even), we have 
\begin{align*}
    \varphi(a_{1}c_{1}, \hdots, a_{m}&c_{m}, a_{m+1}, \hdots, a_{2m})\\
    &=\varphi(a_{1}, \hdots, a_{m}, c_{m}a_{m+1}, \hdots, c_{1}a_{2m}),
\end{align*} 
for all $a_{1}, a_{2}, \hdots, a_{2m}, c_{1}, \hdots, c_{m} \in \mathcal{A}$.
\end{enumerate}
 \end{defn}
 Now, we fix some notations to represent such maps.
 \begin{enumerate}
 \item $ICC_{loc}^{s}(\mathcal{A}^k, \mathcal{B})$:  The space of all symmetric local completely contractive invariant $k$-linear maps from $\mathcal{A}^k$ to $\mathcal{B}.$
 \item $ICP_{loc}^{s}(\mathcal{A}^k,\mathcal{B})$:  The set of all symmetric local completely positive invariant $k$-linear maps from $\mathcal{A}^k$ to $\mathcal{B}.$
 \item $ICPCC_{loc}^{s}(\mathcal{A}^k,\mathcal{B})$:  The set of all symmetric local completely contractive and local completely positive invariant $k$-linear maps from $\mathcal{A}^k$ to $\mathcal{B}.$
  \end{enumerate}
 \subsection{Quantized domain and Representation of local $C^*$-algebras} A quantized domain in a Hilbert space $\mathcal{H}$ is represented by a triple $\{\mathcal{H}; \mathcal{E}; \mathcal{D}\}$, where $\mathcal{E}=\{ \mathcal{H}_{\ell}: \ell\in \Omega\}$ is a directed upward filtered family of closed subspaces such that $\mathcal{D}=\bigcup\limits_{\ell\in \Omega}\mathcal{H}_{\ell}$ is dense in $\mathcal{H}.$ Like Hilbert space, quantized domain enjoys direct sum operations: That is, if $\mathcal{H}^{n}, \mathcal{H}_{\ell}^{n}, \mathcal{D}^{n}$ are the $n$-direct sum copies of $\mathcal{H}, \mathcal{H}_{\ell}, \mathcal{D}$ respectively and 
 $\mathcal{E}^{n}=\{\mathcal{H}^{n}_{\ell}: \ell\in \Omega\}$ for $n\in \N$, then  $\{\mathcal{H}^{n}; \mathcal{E}^{n};\mathcal{D}
 ^{n}\}$ represents quantized domain in the Hilbert space $\mathcal{H}^{n}$.  Let $\{\mathcal{H}^{(i)}; \mathcal{E}^{(i)};\mathcal{D}^{(i)}\}$ be a quantized domain in Hilbert space $\mathcal{H}^{(i)}$ and $\mathcal{E}^{(i)}= \{\mathcal{H}^{(i)}_{\ell}: \ell\in \Omega\}$ for $i=1,2$. Then the notation $\mathcal{E}^{(1)} \subseteq \mathcal{E}^{(2)}$ means that $\mathcal{H}_{\ell}^{(1)} \subseteq \mathcal{H}_{\ell}^{(2)}$ for all $\ell \in \Omega$.  
 
 Let $\mathcal{H, K}$ be two Hilbert spaces. A linear operator $T: {\rm dom}(T)\subseteq \mathcal{H}\to \mathcal{K}$ is called \emph{densely defined} if ${\rm dom}(T)$ is dense in $\mathcal{H}.$ Then adjoint of $T$ is a linear map $T^{\bigstar}: {\rm dom}(T^{\bigstar})\to \mathcal{H}$ with domain ${\rm dom}(T^{\bigstar})$ is given by
  \begin{equation*}
  {\rm dom}(T^{\bigstar})= \big\{\xi\in \mathcal{K}: \eta\mapsto \langle T\eta, \xi \rangle_{\mathcal{K}} \text{ is continuous, for every}\; \eta\in {\rm dom}(T) \big\}.
\end{equation*}
 It implies that
 \begin{equation*}
 \langle T\eta, \xi \rangle_{\mathcal{K}}=\langle \eta, T^{\bigstar}\xi \rangle_{\mathcal{H}},\; \text{ for all }\; \eta\in {\rm dom}(T) \; \text{and}\; \xi\in {\rm dom}(T^{\bigstar}).
 \end{equation*}
 Let  $\{\mathcal{H}; \mathcal{E}; \mathcal{D}\}$ be a quantized domain in the Hilbert space $\mathcal{H}$ where $\mathcal{E}=\{\mathcal{H}_{\ell}: \ell\in \Omega\}$.  Let $P_{\ell}$ be the projection of $\mathcal{H}$ onto $\mathcal{H}_{\ell}$ for each $\ell\in \Omega$ and $\mathcal{L(D)}$ denote the set of all linear operators on $\mathcal{D}.$ The algebra of all \emph{non-commutative continuous functions on a quantized domain} is defined as
 \begin{align}
 C^*_{\mathcal{E}}(\mathcal{D})=\big\{T\in \mathcal{L(D)}: TP_{\ell}= P_{\ell}TP_{\ell}, P_{\ell}T\subseteq TP_{\ell} , \ell\in \Omega\big\}.
 \end{align}
 It is easy to see that $C^*_{\mathcal{E}}(\mathcal{D})$ is a $\ast$-algebra, where the $*$-operation is defined by  $$T^{*}:= T^{\bigstar}|_{\mathcal{D}},\; \text{ for all } T\in C^*_{\mathcal{E}}(\mathcal{D}).$$
 Further, $C^*_{\mathcal{E}}(\mathcal{D})$ is a locally $C^*$-algebra determined by the family   $\{\Vert .\Vert_{\ell} : \ell\in \Omega\}$ of $C^*$-seminorms defined by
 \begin{equation*}
  \Vert T \Vert_{\ell} := \Vert T|_{\mathcal{H}_{\ell}}\Vert \text{~for all~}  T\in C^*_{\mathcal{E}}(\mathcal{D}).
 \end{equation*}
 In 2010, Anar Dosiev obtained GNS type representation  for locally $C^*$-algebras \cite[Theorem 7.2]{Dosiev}. The result states that every locally $C^{\ast}$-algebra $\mathcal{A}$ is identified with a $*$-subalgebra of $C^{\ast}_{\mathcal{E}}(\mathcal{D})$ for some quantaized domain $\{\mathcal{H}; \mathcal{E}; \mathcal{D}\}$.
Also, the author studied the Stinespring type representation for local completely contractive and local completely positive maps $\varphi \colon \mathcal{A}\to C^{*}_{\mathcal{E}}(\mathcal{D}).$ It is regarded as a structure theorem for such maps in the sense that $\varphi$ can be dilated to a unital local contractive $\ast$-homomorphism \cite{Dosiev}. Detailed explanation and motivation for studying these maps can be found in \cite{Dosiev12}.

\section{Stinespring type representation}
In this section, we obtain Stinespring type representation for symmetric local completely contractive and local completely positive invariant $k$-linear maps. Our result is a locally convex analogue of Heo's work  \cite{Heo} as we deal with a multilinear local positive map between local $C^{\ast}$-algebras. In particular, if $k =1$ then the result coincides with \cite[Theorem 5.1]{Dosiev}.

\begin{lemma}\label{Lemma: semi inner product} Let $\varphi \in ICPCC_{loc}^{s}(\mathcal{A}^{k}, {C}^{\ast}_{\mathcal{E}}(\mathcal{D}))$ and $m= [\frac{k+1}{2}]$. Then there exist a sesquilinear form on $\mathcal{A}^{\otimes m}\otimes \mathcal{D}$ and a quantized domain $\{\mathcal{H}^{\varphi}; \mathcal{E}^{\varphi};\mathcal{D}^{\varphi}\}$associated to $\varphi$.
\end{lemma}
\begin{proof} 
Let $\mathcal{A}^{\otimes m}\otimes \mathcal{D}$ be an algebraic tensor product of $\mathcal{A}^{\otimes m}$ with $\mathcal{D}$, where 
$\mathcal{A}^{\otimes m}$  denotes the $m$-fold tensor product of $\mathcal{A}$ with itself. Now we define a map $\big\langle \cdot, \cdot \big\rangle_{\varphi} \colon (\mathcal{A}^{\otimes m}\otimes \mathcal{D}) \times (\mathcal{A}^{\otimes m}\otimes \mathcal{D}) \to \mathbb{C}$ by
\begin{align*}\label{Equation: positive semi definite 1}
\Big\langle &\sum_{j=1}^{n} a_{1_{j}}\otimes \hdots \otimes a_{m_{j}}\otimes g_{j}, \sum_{i=1}^{t} b_{1_{i}}\otimes \hdots \otimes b_{m_{i}}\otimes h_{i} \Big\rangle_{\varphi} \\
&:=\left\{ \begin{array}{cc}
     \sum\limits_{j,i=1}^{n,t}\big\langle \varphi(b_{m_{i}}^*, \hdots, b_{2_{i}}^*, b_{1_{i}}^*a_{1_{j}}, a_{2_{j}}, \hdots, a_{m_{j}})g_{j},\; h_{i}\big\rangle_{\mathcal{D}} &\text{if} \; k = 2m-1 \\
   &\\
    \sum\limits_{j,i=1}^{n,t}\big\langle \varphi(b_{m_{i}}^*, \hdots, b_{2_{i}}^*, b_{1_{i}}^*, a_{1_{j}}, a_{2_{j}}, \hdots, a_{m_{j}})g_{j},\; h_{i}\big\rangle_{\mathcal{D}} &\text{if} \; k = 2m, \\
\end{array} \right.
\end{align*}
for all $a_{p_{j}}, b_{p_{i}} \in \mathcal{A}; g_{j}, h_{i} \in \mathcal{D},\;  j\in\{1, \hdots, n\}, i\in \{1, \hdots, t\}$ and $p \in \{1,\hdots,m\}$. Since $\varphi$ is $k$-linear, it follows that $\langle \cdot , \cdot \rangle_{\varphi}$ is linear in the first variable and  anti-linear in the second variable. Firstly, we show that  $\langle\cdot,\cdot\rangle_{\varphi}$ is a positive semi-definite map. Let  $a_{p_{j}}\in \mathcal{A}, g_{j}\in \mathcal{D}, p\in \{1, \hdots, m\}$ and 
$j \in \{1, \hdots, n\}$.  Since $\mathcal{D}$ is the union space of upward filtered family $\mathcal{E}$, there is an $\ell \in \Omega$ such that $g_{j} \in \mathcal{H}_{\ell}$ for $1 \leq j \leq n.$ Since $\varphi \in ICPCC_{loc}^{s}(\mathcal{A}^{k}, {C}^{\ast}_{\mathcal{E}}(\mathcal{D}))$ and each $\mathcal{H}_{\ell}$ is invariant under $\varphi(b_{1}, b_{2}, \hdots, b_{k})$ for all $b_{1},b_{2}, \hdots, b_{k}\in \mathcal{A},$
  we have
\begin{align}
\notag \Big\langle &\sum_{j=1}^{n} a_{1_{j}}\otimes \hdots \otimes a_{m_{j}}\otimes g_{j}, \sum_{j=1}^{n} a_{1_{j}}\otimes \hdots \otimes a_{m_{j}}\otimes g_{j}\Big \rangle_{\varphi} \\
\notag &=\left\{ \begin{array}{cc}
     \sum\limits_{j,i=1}^{n}\big\langle \varphi(a_{m_{i}}^*, \hdots, a_{2_{i}}^*, a_{1_{i}}^*a_{1_{j}}, a_{2_{j}}, \hdots, a_{m_{j}})g_{j},\; g_{i}\big\rangle_{\mathcal{H}_{\ell}} &\text{if}\; \; k = 2m-1 \\
   &\\
    \sum\limits_{j,i=1}^{n}\big\langle \varphi(a_{m_{i}}^*, \hdots, a_{2_{i}}^*, a_{1_{i}}^*, a_{1_{j}}, a_{2_{j}}, \hdots, a_{m_{j}})g_{j},\; g_{i}\big\rangle_{\mathcal{H}_{\ell}} &\text{if} \;\; k = 2m \\
\end{array} \right.\\
  \notag &=\left\{ \begin{array}{cc}
    \Bigg\langle \Bigg[
     \varphi(a_{m_i}^*, \hdots, a_{2_i}^*, a_{1_i}^*a_{1_j}, a_{2_j}, \hdots, a_{m_j})\Bigg]_{i,j=1}^{n} \begin{bmatrix}g_{1}\\\vdots\\ g_{n} \end{bmatrix}, \begin{bmatrix}g_{1}\\\vdots\\ g_{n} \end{bmatrix}\Bigg\rangle_{\mathcal{H}_{\ell}^{n}} &\text{if}\; \; k = 2m-1 \\
   &\\
    \Bigg\langle \Bigg[
     \varphi(a_{m_i}^*, \hdots, a_{2_i}^*, a_{1_i}^*, a_{1_j}, a_{2_j}, \hdots, a_{m_j})\Bigg]_{i,j=1}^{n} \begin{bmatrix}g_{1}\\\vdots\\ g_{n} \end{bmatrix}, \begin{bmatrix}g_{1}\\\vdots\\ g_{n} \end{bmatrix}\Bigg\rangle_{\mathcal{H}_{\ell}^{n}} &\text{if} \;\; k = 2m \\
\end{array} \right.\\
 \notag &=\left\{ \begin{array}{cc}
    \Big \langle 
     \varphi_{n}(D_{m}^*, \hdots, D_{2}^*, R_{1}^*R_{1}, D_{2}, \hdots, D_{m}) \tilde{g},\; \tilde{g}\Big\rangle_{\mathcal{H}_{\ell}^{n}} &\text{if} \;\; k = 2m-1 \\
   &\\  
       \Big\langle\varphi_{n}(D_{m}^*, \hdots, D_{2}^*, R_{1}^*, R_{1}, D_{2}, \hdots, D_{m}) \tilde{g},\; \tilde{g}\Big\rangle_{\mathcal{H}_{\ell}^{n}} &\text{if} \;\; k = 2m ,
\end{array} \right. 
\end{align}
where
\begin{equation}\label{Eq: Essential}
    \tilde{g}= 
\begin{bmatrix}
g_{1}\\
\vdots\\
g_{n}
\end{bmatrix}, \; D_{p}= \begin{bmatrix}
a_{p_1}& 0 &\hdots &0\\
0& a_{p_2}&\hdots &0\\
\vdots&\vdots&\ddots&\vdots\\
0&0&\hdots &a_{p_n}
\end{bmatrix} \; \text{and}\; R_{1}= 
\begin{bmatrix}
a_{1_1}& a_{1_2}\hdots& a_{1_n}\\
0     &\hdots       &0\\
\vdots&\ddots&\vdots\\
0     &\hdots       &0\\
\end{bmatrix},
\end{equation}  
for $2 \leq p\leq m$. If $k = 2m-1$, then 
\begin{equation*}
    (D_{m}^*, \hdots, D_{2}^*, R_{1}^*R_{1}, D_{2}, \hdots, D_{m})^{\ast} = (D_{m}^*, \hdots, D_{2}^*, R_{1}^*R_{1}, D_{2}, \hdots, D_{m}) \; \text{and}\; R_{1}^{\ast}R_{1} \geq_{\alpha} 0,
\end{equation*}  
for every $\alpha \in \Lambda$. If $k = 2m$, then 
\begin{equation*}
    (D_{m}^*, \hdots, D_{2}^*, R_{1}^*, R_{1}, D_{2}, \hdots, D_{m})^{\ast} = (D_{m}^*, \hdots, D_{2}^*, R_{1}^*, R_{1}, D_{2}, \hdots, D_{m}).
\end{equation*}
Since $\varphi$ is a local completely positive map, it follows that
\begin{equation*}
    \Big\langle \sum_{j=1}^{n} a_{1_{j}}\otimes \hdots \otimes a_{m_{j}}\otimes g_{j}, \sum_{j=1}^{n} a_{1_{j}}\otimes \hdots \otimes a_{m_{j}}\otimes g_{j}\Big \rangle_{\varphi} \geq 0.
\end{equation*}
   Now by using the symmetric nature of $\varphi,$ we show that $\big\langle \cdot, \cdot \big\rangle_{\varphi}$ satisfies Hermitian property. Let $a_{p_{j}}, b_{p_{i}} \in \mathcal{A}, \; g_{j}, h_{i} \in \mathcal{D}$ for $j\in \{1,\hdots, n\},\; i \in \{1, \hdots, t\}$ and $p\in \{1, \hdots, m\}.$ Then 
   \begin{align*}
&\Big \langle\sum_{j=1}^n a_{1_j}\otimes\hdots\otimes a_{m_j}\otimes g_{j},\; \sum_{i=1}^t b_{1_i}\otimes\hdots \otimes b_{m_i}\otimes h_{i} \Big\rangle_{\varphi}\\
&=\left\{ \begin{array}{cc}
     \sum\limits_{j,i=1}^{n,t} \Big\langle\varphi(b_{m_i}^*,\hdots, b_{2_i}^{\ast}, b_{1_i}^*a_{1_j}, a_{2_j}, \hdots, a_{m_j})g_{j},\; h_{i} \Big\rangle_{\mathcal{H}_{\ell}}  &\text{if} \; k = 2m-1 \\
     \sum\limits_{j,i=1}^{n,t} \Big\langle\varphi(b_{m_i}^*,\hdots, b_{2_i}^{\ast}, b_{1_i}^*,a_{1_j}, a_{2_j}, \hdots, a_{m_j})g_{j},\; h_{i} \Big\rangle_{\mathcal{H}_{\ell}} &\text{if} \; k = 2m \\
\end{array} \right.\\
&=\left\{ \begin{array}{cc}
     \sum\limits_{j,i=1}^{n,t} \Big\langle g_{j}, \varphi(b_{m_i}^*,\hdots, b_{2_i}^{\ast}, b_{1_i}^*a_{1_j}, a_{2_j}, \hdots a_{m_j})^*h_{i} \Big\rangle_{\mathcal{H}_{\ell}} &\text{if} \; k = 2m-1 \\
     \sum\limits_{j,i=1}^{n,t} \Big\langle g_{j}, \varphi(b_{m_i}^*,\hdots, b_{2_i}^{\ast}, b_{1_i}^*, a_{1_j}, a_{2_j}, \hdots a_{m_j})^*h_{i} \Big\rangle_{\mathcal{H}_{\ell}} &\text{if} \; k = 2m \\
\end{array} \right.\\
&=\left\{ \begin{array}{cc}
     \sum\limits_{j,i=1}^{n,t} \Big\langle g_{j},\; \varphi(a_{m_{j}}^{\ast}, \hdots, a_{2_{j}}^{\ast}, a_{1_{j}}^{\ast}b_{1_{i}}, b_{2_{i}}, \hdots, b_{m_i})h_{i} \Big\rangle_{\mathcal{H}_{\ell}} &\text{if} \; k = 2m-1 \\
     \sum\limits_{j,i=1}^{n,t} \Big\langle g_{j},\; \varphi(a_{m_{j}}^{\ast}, \hdots, a_{2_{j}}^{\ast}, a_{1_{j}}^{\ast}, b_{1_{i}}, b_{2_{i}}, \hdots, b_{m_i})h_{i} \Big\rangle_{\mathcal{H}_{\ell}} &\text{if} \; k = 2m \\
\end{array} \right.\\
&=\left\{ \begin{array}{cc}
     \sum\limits_{j,i=1}^{n,t} \overline{\Big\langle  \varphi(a_{m_{j}}^{\ast}, \hdots, a_{2_{j}}^{\ast}, a_{1_{j}}^{\ast}b_{1_{i}}, b_{2_{i}}, \hdots, b_{m_i})h_{i}, \; g_{j} \Big\rangle}_{\mathcal{H}_{\ell}} &\text{if} \; k = 2m-1 \\
     \sum\limits_{j,i=1}^{n,t} \overline{\Big\langle \varphi(a_{m_{j}}^{\ast}, \hdots, a_{2_{j}}^{\ast}, a_{1_{j}}^{\ast}, b_{1_{i}}, b_{2_{i}}, \hdots, b_{m_i})h_{i}, \; g_{j} \Big\rangle}_{\mathcal{H}_{\ell}} &\text{if} \; k = 2m \\
\end{array} \right.\\
&=\overline{\Big \langle\sum_{i=1}^{t} b_{1_i}\otimes\hdots \otimes b_{m_i}\otimes h_{i},\sum_{j=1}^{n} a_{1_j}\otimes\hdots\otimes a_{m_j}\otimes g_{j} \Big\rangle_{\varphi}}.
\end{align*}
Thus $\langle \cdot,\cdot\rangle_{\varphi}$ is a semi-inner product on $\mathcal{A}^{\otimes m}\otimes \mathcal{D}$. To prove the second part of the claim, let us define $\mathcal{N}_{\varphi} := \big\{\xi \in \mathcal{A}^{\otimes m} \otimes \mathcal{D}:\; \langle \xi, \xi \rangle_{\varphi} = 0\big\}.$ Then by the Cauchy-Schwartz inequality of semi-inner prouct $\langle \cdot, \cdot \rangle_{\varphi}$ on $\mathcal{A}^{\otimes m} \otimes \mathcal{D}$, we have $ \mathcal{N}_{\varphi} = \big\{\xi \in \mathcal{A}^{\otimes m}\otimes \mathcal{D}:\; \big\langle \xi, \eta \big\rangle_{\varphi} = 0, \; \text{for every} \; \eta \in \mathcal{A}^{\otimes m}\otimes \mathcal{D}\big\}.$
This shows that $\mathcal{N}_{\varphi}$ is a closed subspace of $\mathcal{A}^{\otimes m}\otimes \mathcal{D}.$ The map $\langle \cdot, \cdot \rangle $ defined by 
\begin{equation*}
\langle \xi + \mathcal{N}_{\varphi},\; \xi^{\prime}+\mathcal{N}_{\varphi}\rangle:= \langle \xi,\; \xi^{\prime}\rangle_{\varphi}, \; \text{for every}\; \xi, \xi^{\prime} \in \mathcal{A}^{\otimes m} \otimes \mathcal{D}
\end{equation*}
is a sesquilinear form on $\mathcal{A}^{\otimes m} \otimes \mathcal{D}.$ Let $\mathcal{H}^{\varphi}$ be the Hilbert space obtained by the completion of $\big(\mathcal{A}^{\otimes m} \otimes \mathcal{D}\big)/ \mathcal{N}_{\varphi}$.  Similarly, for each $\ell \in \Omega$ we define the Hilbert space $\mathcal{H}_{\ell}^{\varphi}$ as,
\begin{equation}\label{Eq: Hlphi}
\mathcal{H}^{\varphi}_{\ell}=\overline{{\rm span}}\{a_{1}\otimes \hdots \otimes a_{m}\otimes g+ \mathcal{N}_{\varphi}: a_{p}\in \mathcal{A}, 1\leq p\leq m, g\in \mathcal{H}_{\ell}\}.
\end{equation}
Since $\mathcal{E}$ is an upward filtered family, it follows from (\ref{Eq: Hlphi})  that $\mathcal{E}^{\varphi}= \{\mathcal{H}_{\ell}^{\varphi}:
\ell \in \Omega\}$ is an upward filtered family of closed subspaces of $\mathcal{H}^{\varphi}$. It remains to show that $\mathcal{D}^{\varphi}=\bigcup\limits_{\ell\in \Omega}\mathcal{H}^{\varphi}_{\ell}$ is dense in $\mathcal{H}^{\varphi}$. Clearly, $\overline{\mathcal{D}^{\varphi}} \subseteq \mathcal{H}^{\varphi}$. Suppose there is a non-zero $\eta\in (\overline{\mathcal{D}^{\varphi}})^{\perp},$ then there exists an $n_{0}\in \mathbb{N}$ such that $\frac{1}{n_{0}} < \Vert \eta\Vert^{2}.$ Since $(\overline{\mathcal{A}^{\otimes m}\otimes
\mathcal{D})/\mathcal{N}_{\varphi}}=\mathcal{H}^{\varphi}$, there exists a $\gamma \in \mathcal{A}^{\otimes
m}\otimes \mathcal{D}/\mathcal{N}_{\varphi}$ of the form $\gamma = \sum\limits_{j=1}^{t}a_{1j}\otimes \hdots \otimes a_{mj}\otimes g_{j}+ \mathcal{N}_{\varphi}$ such that $\big\| \eta- \gamma \big\|^{2} < \frac{1}{n_{0}}$. Without loss of generality, for every $1\leq j \leq t$ we have $g_{j}\in \mathcal{H}_{\ell}$ for some $\ell \in \Omega$. It follows that $\gamma \in  \mathcal{H}_{\ell}^{\varphi} \subseteq \mathcal{D}^{\varphi}$ and hence $\langle \eta, \gamma \rangle=0.$ Therefore, 
\begin{equation*}
\Vert \eta\Vert^{2}\leq \Vert \eta\Vert^{2}+\|\gamma\|^2 = \|\eta - \gamma\|^2 <\frac{1}{n_{0}}.
\end{equation*}
This contradicts the fact that $\frac{1}{n_{0}}< \Vert \eta\Vert^{2}.$ Thus $\mathcal{D}^{\varphi}$ is dense in $\mathcal{H}^{\varphi}$. As a result, $\{\mathcal{H}^{\varphi};\mathcal{E}^{\varphi}; \mathcal{D}^{\varphi}\}$ is a quantized domain associated to $\varphi$.
\end{proof}
Next, we obtain Stinespring type representation for  symmetric local completely contractive and local completely positive invariant $k$-linear maps. The original construction of the Stinespring theorem for CP-maps on $C^*$-algebras can be found in \cite{Stinespring}.
\begin{thm}\label{Theorem: Multilinear Stinespring}
Let $\varphi \in ICPCC_{loc}^{s}(\mathcal{A}^{k}, {C}^{\ast}_{\mathcal{E}}(\mathcal{D}))$. Then there is a quantized domain $\big\{\mathcal{H}^{\varphi}; \mathcal{E}^{\varphi}; \mathcal{D}^{\varphi}\big\}$, a contraction $V_{\varphi} \in \mathcal{B}(\mathcal{H}, \mathcal{H}^{\varphi})$ and a commuting family $\{\pi_{p}^{\varphi}: 1 \leq p \leq m\}$ of unital local contractive $\ast$-homomorphisms from $\mathcal{A}$ to $C^{\ast}_{\mathcal{E}^{\varphi}}(\mathcal{D}^{\varphi})$ such that $V_{\varphi}(\mathcal{E}) \subseteq \mathcal{E}^{\varphi}$ and 
\begin{enumerate}
    \item[(i)] if $k = 2m-1$ odd, then 
    \begin{equation*}
        \varphi(a_{1}, \hdots, a_{m-1}, a_{m}, a_{m+1}, \hdots, a_{2m-1}) \subseteq V_{\varphi}^{\ast}\pi^{\varphi}_{1}(a_{m})\pi^{\varphi}_{2}(a_{m-1}a_{m+1}) \cdots \pi^{\varphi}_{m}(a_{1}a_{2m-1})V_{\varphi},
    \end{equation*}
    for all $a_{1}, a_{2}, \hdots, a_{2m-1} \in \mathcal{A}$.
    \item[(ii)] if $k = 2m$ even, then 
    \begin{equation*}
        \varphi(a_{1}, \hdots, a_{m-1}, a_{m}, a_{m+1}, \hdots, a_{2m}) \subseteq V_{\varphi}^{\ast}\pi^{\varphi}_{1}(a_{m}a_{m+1})\pi^{\varphi}_{2}(a_{m-1}a_{m+2}) \cdots \pi^{\varphi}_{m}(a_{1}a_{2m})V_{\varphi},
    \end{equation*}
    for all $a_{1}, a_{2}, \hdots, a_{2m} \in \mathcal{A}$.
\end{enumerate}
\end{thm}
\begin{proof} By Lemma \ref{Lemma: semi inner product}, there is a quantized domain $\{\mathcal{H}^{\varphi}; \mathcal{E}^{\varphi}; \mathcal{D}^{\varphi}\}$ associated to $\varphi$. Recall that $\mathcal{H}^{\varphi}$ is obtained by the completion of the quotient space $(\mathcal{A}^{\otimes m} \otimes \mathcal{D})/\mathcal{N}_{\varphi}$. For $1\leq p\leq m$ and  $a\in \mathcal{A},$ we define a linear map
$\pi^{\varphi}_{p}(a):(\mathcal{A}^{\otimes m}\otimes \mathcal{D})/\mathcal{N}_{\varphi} \to (\mathcal{A}^{\otimes m}\otimes \mathcal{D})/\mathcal{N}_{\varphi}$ by
\begin{align}
  \notag \pi_{p}^{\varphi}(a)\Big(\sum_{j=1}^n a_{1_{j}}\otimes\hdots\otimes a_{p_{j}}\otimes & \hdots\otimes a_{m_{j}}\otimes g_{j}+\mathcal{N}_{\varphi} \Big)\\
 &=\sum_{j=1}^n a_{1_{j}}\otimes\hdots \otimes aa_{p_{j}}\otimes \hdots\otimes a_{m_{j}}\otimes g_{j}+\mathcal{N}_{\varphi}
  \end{align}
  for all $a_{p_{j}} \in \mathcal{A},\; g_{j} \in \mathcal{D},\; j \in \{1,\hdots, n\}.$ Now we show that $\pi^{\varphi}_{p}$ is a well defined linear map. We can choose $\ell \in \Omega$ with $g_{j} \in \mathcal{H}_{\ell}$ for $j \in \{1,\hdots, n\}$. Then
  
\begin{align} \label{Eq: order relation on phi}
 \notag  &\Big\| \pi_{p}^{\varphi}(a)\big(\sum_{j=1}^n a_{1_j}\otimes \hdots \otimes a_{p_j}\otimes \hdots \otimes a_{m_j}\otimes g_{j}+\mathcal{N}_{\varphi}\big)\Big\|_{\mathcal{H}_{\ell}^{\varphi}}^2\\
 \notag  &=\Big\| \sum_{i=1}^n a_{1_{j}}\otimes \hdots \otimes aa_{p_{j}}\otimes \hdots\otimes a_{m_{j}}\otimes g_{j}+\mathcal{N}_{\varphi}\Big\|_{\mathcal{H}_{\ell}^{\varphi}}^2\\
   \notag  &=\left\{ \begin{array}{cc}
   \sum\limits_{i,j=1}^n \big \langle \varphi(a_{m_i}^*, \hdots,  a_{p_i}^*a^*, \hdots, a_{2_i}^*,a_{1_i}^*a_{1_j}, \hdots ,aa_{p_j}, \hdots, a_{m_j})g_{j}, g_{i}\big \rangle_{\mathcal{H}_{\ell}}  &\text{if} \; k = 2m-1 \\
   \sum\limits_{i,j=1}^n \big \langle \varphi(a_{m_i}^*, \hdots,  a_{p_i}^*a^*, \hdots, a_{2_i}^*,a_{1_i}^*, a_{1_j}, \hdots ,aa_{p_j}, \hdots, a_{m_j})g_{j}, g_{i}\big \rangle_{\mathcal{H}_{\ell}}&\text{if} \; k = 2m \\
\end{array} \right.\\
  &=\left\{ \begin{array}{cc}
   \big\langle \varphi_{n}(D_{m}^*, \hdots, D_{p}^* (a^* \otimes I_{n}), \hdots, R_{1}^*R_{1}, \hdots ,(a\otimes I_{n}) D_{p}, \hdots, D_{m})\tilde{g}, \tilde{g}\big\rangle_{\mathcal{H}^{n}_{\ell}}  &\text{if} \; k = 2m-1 \\
   &\\
 \big \langle \varphi_{n}(D_{m}^*, \hdots, D_{p}^* (a^* \otimes I_{n}), \hdots,  R_{1}^*, R_{1}, \hdots ,(a\otimes I_{n}) D_{p}, \hdots, D_{m})\tilde{g}, \tilde{g}\big\rangle_{\mathcal{H}^{n}_{\ell}} &\text{if} \; k = 2m, \\
\end{array} \right.
\end{align}
where $\widetilde{g}, R_{1}$ and $D_{p}$ for $2 \leq p \leq m$ are defined as in Equation \eqref{Eq: Essential}. 
Since $\varphi$ is a local completely positive map, for given $\ell\in \Omega$, there exists an $\alpha_{\ell}\in \Lambda$ such that $  \varphi_{n}(B_{1}, B_{2},\hdots , B_{k})\geq_{\ell} 0 $ whenever $(B_{1}, B_{2},\hdots ,B_{k}) =_{\alpha_{\ell}}(B_{k}^*,\hdots, B_{2}^*, B_{1}^*)$ and $\varphi_{n}(B_{1}, B_{2}, \hdots, B_{k})=_{\ell}0 $ whenever $B_{i}=_{\alpha_{\ell}}0$  for some $1\leq i\leq k,$ and additionally
$B_{m}\geq_{\alpha_{\ell}}0$ if $k=2m-1$ is odd.
 Let us consider the case $k = 2m-1$. Notice that $\pi_{\alpha_{\ell}}(p_{\alpha_{\ell}}(a)^2I-a^*a)$ is positive in  $\mathcal{A}_{\alpha_{\ell}},$ thus there exists an element $b\in \mathcal{A}$ such that $\pi_{\alpha_{\ell}}(p_{\alpha}(a)^2I-a^*a)=\pi_{\alpha_{\ell}}(b^*b).$ Therefore $p_{\alpha_{\ell}}(a)^2I-a^*a= b^*b+ c$ for some $c\in \mathcal{N}_{\alpha_{\ell}}.$
As $c=_{\alpha_{\ell}}0,$ we have $(c\otimes I_{n})D_{p}=_{\alpha_{\ell}} 0.$ Applying these facts on $\varphi_n,$ notice that
  $\varphi_{n}(D_{m}^*, \hdots, D_{p}^*, \hdots, D_{2}^*, R_{1}^*R_{1},
 D_{2},\hdots , (c \otimes I_{n}) D_{p}, \hdots, D_{m})=_{\ell}0.$ 
  Applying this argument as well as the invariant property of $\varphi$ carefully, we obtain
\begin{align*}
&\Big\langle \varphi_{n}(D_{m}^*, \hdots, D_{p}^*, \hdots, D_{2}^*, R_{1}^*R_{1},
 D_{2},\hdots , ((p_{\alpha_{\ell}}(a)^2I- a^*a) \otimes I_{n}) D_{p}, \hdots, D_{m})\tilde{g}, \; \tilde{g}\Big\rangle_{\mathcal{H}^{n}_{\ell}}\\
 =&\Big\langle \varphi_{n}(D_{m}^*, \hdots, D_{p}^*, \hdots, D_{2}^*, R_{1}^*R_{1},
 D_{2},\hdots , (b^*b \otimes I_{n}) D_{p}, \hdots, D_{m})\tilde{g}, \; \tilde{g}\Big\rangle_{\mathcal{H}^{n}_{\ell}}\\
=&\Big\langle \varphi_{n}(D_{m}^*, \hdots, D_{p}^*(b^* \otimes I_{n}), \hdots, D_{2}^*, R_{1}^*R_{1},
  D_{2},\hdots ,(b \otimes I_{n}) D_{p}, 
   D_{m})\tilde{g},\; \tilde{g}\Big\rangle_{\mathcal{H}^{n}_{\ell}}\\
  &\geq 0.
\end{align*}
This implies that 
 \begin{align*}
 \Big\langle \varphi_{n}(D_{m}^*, \hdots,& D_{p}^*, \ldots, D_{2}^*, R_{1}^*R_{1},
 D_{2},\hdots , (a^*a \otimes I_{n}) D_{p}, \hdots, D_{m})\tilde{g},\; \tilde{g}\Big\rangle_{\mathcal{H}^{n}_{\ell}} \\
 &\leq   p_{\alpha_{\ell}}(a)^2\Big\langle \varphi_{n}(D_{m}^*, \hdots, D_{p}^*, \hdots, D_{2}^*, R_{1}^*R_{1},
 D_{2},\hdots ,  D_{p}, \hdots, D_{m})\tilde{g},\;  \tilde{g}\Big\rangle_{\mathcal{H}^{n}_{\ell}}\\
 & = p_{\alpha_{\ell}}(a)^2 \sum\limits_{i,j=1}^n  \big\langle \varphi(a_{m_i}^*, \hdots, a_{2_i}^*,a_{1_i}^*a_{1_j}, a_{2_{j}}, \hdots , a_{m_j})g_{j},\; g_{i}\big\rangle_{\mathcal{H}_{\ell}}\\
 &= p_{\alpha_{\ell}}(a)^2 \Big\|\sum\limits_{j=1}^{n} a_{1_{j}} \otimes a_{2_{j}} \otimes \ldots \hdots \otimes a_{m_{j}}\otimes g_{j} + \mathcal{N}_{\varphi}\Big\|_{\mathcal{H}_{\ell}^{\varphi}}^{2}
 \end{align*}
 Similarly, the above inequality holds true for $k = 2m.$ From Equation \eqref{Eq: order relation on phi}, we conclude that 
 \begin{equation}\label{Eq: pi's are local contractive}
 \|\pi_{p}^{\varphi}(a)\|_{\mathcal{H}_{\ell}^{\varphi}} \leq p_{\alpha_{\ell}}(a); \text{ for every}\; a \in \mathcal{A},\; \ell \in \Omega.
 \end{equation}
 Equivalently, $\pi_{p}^{\varphi}: \mathcal{A} \to C_{\mathcal{E}^{\varphi}}^*(\mathcal{D}^{\varphi})$ is a local contractive map. Next we show that $a \mapsto \pi_{p}^{\varphi}(a)$ is a unital $\ast$-homomorphism. Let $a_{p_{j}}, b_{p_{i}} \in \mathcal{A} $ and $g_{j}, h_{i} \in \mathcal{H}_{\ell}$ for $p \in \{1, \hdots, m\}, j\in \{1, \hdots, t\}, i \in \{1, \hdots, n\}$. Then by using invariant property of $\varphi$, we get\\
 \begin{align*}
 \notag &\Big\langle \pi_{p}^{\varphi}(a)^{\ast}\Big(\sum_{j=1}^t a_{1_j}\otimes \hdots \otimes a_{p_j}\otimes\hdots \otimes a_{m_j}\otimes g_{j}+\mathcal{N}_{\varphi}\Big),\;  \sum_{i=1}^{n}b_{1_i} \otimes\hdots\otimes b_{p_i}\otimes\hdots\otimes b_{m_i}\otimes h_{i}+\mathcal{N}_{\varphi} \Big\rangle_{\mathcal{H}_{\ell}^{\varphi}}\\ 
 \notag &=\Big\langle \sum_{j=1}^t a_{1_j}\otimes \hdots \otimes a_{p_j}\otimes\hdots \otimes a_{m_j}\otimes g_{j}+\mathcal{N}_{\varphi},\; \pi_{p}^{\varphi}(a) \Big(\sum_{i=1}^{n}b_{1_i} \otimes\hdots\otimes b_{p_i} \otimes \hdots\otimes b_{m_i}\otimes h_{i}+\mathcal{N}_{\varphi}\Big) \Big\rangle_{\mathcal{H}_{\ell}^{\varphi}}\\
\notag &= \left\{ \begin{array}{cc}
          \sum\limits_{i,j=1}^{n, t}\Big\langle \varphi(b_{m_i}^*, \hdots, b_{p_i}^*a^{\ast}, \hdots, b_{1_i}^*a_{1_j}, \hdots,  a_{p_j}, \hdots, a_{m_j}) g_{j} , h_{i}\Big \rangle_{\mathcal{H}_{\ell}} &\; \text{if}\; k = 2m-1  \\
      \notag      \sum\limits_{i,j=1}^{n, t}\Big\langle \varphi(b_{m_i}^*, \hdots, b_{p_i}^* a^{\ast}, \hdots, b_{1_i}^*,a_{1_j}, \hdots,  a_{p_j}, \hdots, a_{m_j}) g_{j} , h_{i}\Big \rangle_{\mathcal{H}_{\ell}}&\; \text{if}\; k = 2m 
    \end{array}\right.\\
\notag &= \left\{ \begin{array}{cc}
          \sum\limits_{i,j=1}^{n, t}\Big\langle \varphi(b_{m_i}^*, \hdots, b_{p_i}^*, \hdots, b_{1_i}^*a_{1_j}, \hdots,  a^{\ast}a_{p_j}, \hdots, a_{m_j}) g_{j} , h_{i}\Big \rangle_{\mathcal{H}_{\ell}} &\; \text{if}\; k = 2m-1  \\
      \notag      \sum\limits_{i,j=1}^{n, t}\Big\langle \varphi(b_{m_i}^*, \hdots, b_{p_i}^* , \hdots, b_{1_i}^*,a_{1_j}, \hdots,  a^{\ast}a_{p_j}, \hdots, a_{m_j}) g_{j} , h_{i}\Big \rangle_{\mathcal{H}_{\ell}}&\; \text{if}\; k = 2m 
    \end{array}\right.\\    
    \notag &=\Big\langle \pi_{p}^{\varphi}(a^{\ast})\Big(\sum_{j=1}^t a_{1_j}\otimes \hdots \otimes a_{p_j}\otimes\hdots \otimes a_{m_j}\otimes g_{j}+\mathcal{N}_{\varphi}\Big),\;  \sum_{i=1}^{n}b_{1_i} \otimes\hdots\otimes b_{p_i}\otimes \hdots\otimes b_{m_i}\otimes h_{i}+\mathcal{N}_{\varphi} \Big\rangle_{\mathcal{H}_{\ell}^{\varphi}}
    \end{align*}
 It follows that  $\pi_{p}^{\varphi}(a)^*|_{\mathcal{H}^{\varphi}_{\ell}}=\pi^{\varphi}_{p}(a^*)|_{\mathcal{H}^{\varphi}_{\ell}}$, for every $\ell \in \Omega$ and hence   $\pi^{\varphi}_{p}(a)^*= \pi_{p}^{\varphi}(a^*)$ for all $a\in \mathcal{A}.$ It implies that $\pi_{p}: \mathcal{A}\to C^*_{\mathcal{E}^{\varphi}}(\mathcal{D}^{\varphi})$ is a $*$-homomorphism and it is local contractive by   Equation \eqref{Eq: pi's are local contractive}. 
 
  Now it remains to show that $\big\{\pi_{p}^{\varphi}\big\}_ {p=1}^{m}$ is a commuting family.   Let $1\leq  p < q \leq m.$ Then  
 \begin{align*}
 \pi_{p}^{\varphi}(a)\pi_{q}^{\varphi}(b)\Big(\sum_{j=1}^{t} a_{1_j}\otimes  \hdots \otimes a_{m_j}&\otimes g_{j}+ \mathcal{N}_{\varphi}\Big)\\
 &=\sum_{j=1}^{t} a_{1_j}\otimes \hdots\otimes aa_{p_j}\otimes \hdots \otimes ba_{q_j}\otimes \hdots\otimes a_{m_j}\otimes g_{j}+ \mathcal{N}_{\varphi}\\
 &=\pi_{q}^{\varphi}(b)\pi_{p}^{\varphi}(a)\Big(\sum_{j=1}^{t} a_{1_j}\otimes \hdots \otimes a_{m_j}\otimes g_{j}+ \mathcal{N}_{\varphi} \Big)
 \end{align*}
 Hence by continuity, we have $\pi_{p}^{\varphi}(a)\pi_{p}^{\varphi}(b)|_{\mathcal{H}_{\ell}^{\varphi}}=\pi_{p}^{\varphi}(b)\pi_{p}^{\varphi}(a)|_{\mathcal{H}_{\ell}^{\varphi}}$ as both the maps agree on the dense subspace of $\mathcal{H}_{\ell}^{\varphi}.$ Since this is true for all $\ell\in \Omega$, we conclude that
 $\pi_{p}^{\varphi}(a)\pi_{p}^{\varphi}(b)=\pi_{p}^{\varphi}(b)\pi_{p}^{\varphi}(a).$

Now we define a map $V_{\varphi}: \mathcal{H} \to \mathcal{H}^{\varphi}$ by $V_{\varphi}(g)= (1\otimes \hdots \otimes 1 \otimes g)+\mathcal{N}_{\varphi}$ \text{~for all~} $g\in \mathcal{H}$. Note that $V_{\varphi}(\mathcal{H}_{\ell})\subseteq \mathcal{H}_{\ell}^{\varphi}$ for each $\ell\in \Omega.$ That is, $V_{\varphi}(\mathcal{E}) \subseteq \mathcal{E}^{\varphi}$. Since $\varphi$ is a local contractive map, for $g \in \mathcal{H}_{\ell}$, we see that $V_{\varphi}$ is a contractive as follows:
\begin{align}\label{Eq: Norm of Vvarphi}
\notag\Vert V_{\varphi}g\Vert_{\mathcal{H}_{\ell}^{\varphi}}^2&= \big\langle (1\otimes\hdots \otimes 1 \otimes g)+ \mathcal{N}_{\varphi},\; (1\otimes\hdots \otimes 1 \otimes g)+ \mathcal{N}_{\varphi} \big\rangle_{\mathcal{H}^{\varphi}}\\
\notag &= \langle \varphi(1, \hdots, 1)g,\; g\rangle_{\mathcal{H}_{\ell}}\\
\notag &\leq \Vert \varphi(1, \hdots, 1)g \Vert_{\mathcal{H}_{\ell}}  \Vert g\Vert_{\mathcal{H}_{\ell}}\\
\notag &\leq \Vert g\Vert^2_{\mathcal{H}_{\ell}}.
\end{align}
So $V_{\varphi}$ is a contractive map on the dense subspace $\mathcal{D}^{\varphi}$ and hence  it is contractive on $\mathcal{H}^{\varphi}$. Thus $V_{\varphi}\colon \mathcal{H}\to \mathcal{H}^{\varphi}$ is a well defined contractive map. Moreover $V_{\varphi}$ is an isometry when $\varphi$ is unital. Suppose $k = 2m-1$, $a_{1}, a_{2}, \hdots, a_{2m-1}\in \mathcal{A}$ and $g, h \in \mathcal{H}_{\ell}$ for some $\ell \in \Omega$, we have  
\begin{align*}
    \Big\langle V_{\varphi}^{\ast}\pi_{1}^{\varphi}(a_m) \pi_{2}^{\varphi}&(a_{m-1}a_{m+1}) \cdots \pi_{m}^{\varphi}(a_{1}a_{2m-1}) V_{\varphi} h,\; g\Big\rangle_{\mathcal{H}_{\ell}}\\
    & = \Big\langle \pi_{1}^{\varphi}(a_m) \pi_{2}^{\varphi}(a_{m-1}a_{m+1}) \cdots \pi_{m}^{\varphi}(a_{1}a_{2m-1}) V_{\varphi} h,\; V_{\varphi}g\Big\rangle_{\mathcal{H}_{\ell}}\\
    &= \Big\langle (a_{m}\otimes a_{m-1}a_{m+1}\otimes \hdots \otimes a_{1}a_{2m-1}\otimes h )+ \mathcal{N}_{\varphi}, \; (1\otimes \hdots \otimes 1 \otimes g)+ \mathcal{N}_{\varphi} \Big\rangle_{\mathcal{H}_{\ell}^{\varphi}}\\
    &= \big\langle \varphi(1, 1, \hdots, a_{m}, a_{m-1}a_{m+1}, \hdots, a_{1}a_{2m-1}) h, \; g \big\rangle_{\mathcal{H}_{\ell}}\\
    &= \big\langle \varphi(a_{1}, a_{2}, \hdots, a_{m},  \hdots, a_{2m-1}) h, \; g\big\rangle_{\mathcal{H}_{\ell}}
\end{align*}
since $\varphi$ is invariant map. Therefore,
\begin{equation*}
        \varphi(a_{1}, \hdots, a_{m-1}, a_{m}, a_{m+1}, \hdots, a_{2m-1}) \subseteq V_{\varphi}^{\ast}\pi^{\varphi}_{1}(a_{m})\pi^{\varphi}_{2}(a_{m-1}a_{m+1}) \cdots \pi^{\varphi}_{m}(a_{1}a_{2m-1})V_{\varphi},
    \end{equation*}
    for all $a_{1}, \hdots, a_{2m-1}\in \mathcal{A}$. Similarly, the result holds for the case $k = 2m$. 
    \end{proof}
\begin{note}
For a given $\varphi \in ICPCC_{loc}^{s}(\mathcal{A}^{k}, {C}^{\ast}_{\mathcal{E}}(\mathcal{D}))$, we denote a Stinspring's triple associated to $\varphi$ by $\big(\{\pi_{p}^{\varphi}\}_{p=1}^m, V_{\varphi}, \{\mathcal{H}^{\varphi}; \mathcal{E}^{\varphi}; \mathcal{D}^{\varphi}\}\big)$, where the terms are defined as in Theorem \ref{Theorem: Multilinear Stinespring}.
\end{note}
\section{Minimality condition}
In this section, we identify a suitable notion of minimality for Stinespring's type representation so as to ensure the uniqueness up to unitary equivalence. In fact,  W. Arveson was the first mathematician to construct a Radon-Nikod\'ym  derivative between two CP-maps \cite{Arv69} by using minimal Stinespring's representation. Here, we discuss the minimality condition of the Stinespring theorem for invariant multilinear local CP-maps and apply this to prove Radon-Nikod\'ym theorem in the next section. Before that we fix some notations for simplicity:
$ \big[ {M}\big] : = \overline{\textrm{span}}\;{M}, \; \text{for any subset}\;M\;$ of a Hilbert space;
\begin{align*} 
\prod\limits_{p=1}^{m}{\pi}_{p}^{\varphi}(a_{{p}_{j}})&:=\pi_{1}^{\varphi}(a_{1_{j}})\pi_{2}^{\varphi}(a_{2_{j}})\hdots  \pi_{m}^{\varphi}(a_{m_{j}}); \text{ and }\\ 
 (\prod\limits_{p=1}^m \pi_{p}^{\varphi}(\mathcal{A})\big)V_{\varphi}\mathcal{D}&:= \Big\{\prod\limits_{p=1}^{m}{\pi}_{p}^{\varphi}(a_{{p}_{j}})V_{\varphi}h_{j}:\;  
  a_{1_{j}}, \hdots, a_{m_{j}} \in \mathcal{A}, h_{j}\in \mathcal{D} \Big\}.\end{align*}
   The following lemma is a key observation to estabalish the minimality condition.
\begin{lemma}\label{Lemma:  Technical} Let  $\varphi \in  ICPCC^{s}_{loc}(\mathcal{A}^k, C^*_{\mathcal{E}}(\mathcal{D}))$ and let $\big(\{\pi_{p}^{\varphi}\}_{p=1}^m, V_{\varphi}, \{\mathcal{H}^{\varphi}; \mathcal{E}^{\varphi}; \mathcal{D}^{\varphi}\}\big)$ be a Stinespring triple 
of the map $\varphi.$ Then
\begin{align}
   \overline{\bigcup_{\ell\in \Omega} \big[(\prod\limits_{p=1}^m \pi^{\varphi}_{p}(\mathcal{A}))V_{\varphi}\mathcal{H}_{\ell}\big]}= \big[(\prod\limits_{p=1}^m
   \pi^{\varphi}_{p}(\mathcal{A}))V_{\varphi}\mathcal{D}\big]
    \end{align}
\end{lemma}
\begin{proof} The proof follows immediately from the fact that $\bigcup\limits_{\ell\in \Omega}\mathcal{H}_{\ell}=\mathcal{D}.$
   \end{proof}
   \begin{prop} 
    Let  $\varphi \in  ICPCC^{s}_{loc}(\mathcal{A}^k, C^*_{\mathcal{E}}(\mathcal{D}))$ and let $\big(\{\pi_{p}^{\varphi}\}_{p=1}^m, V_{\varphi}, \{\mathcal{H}^{\varphi}; \mathcal{E}^{\varphi}; \mathcal{D}^{\varphi}\}\big)$ be a Stinespring triple 
of $\varphi.$  Then one can construct another Stinespring triple $\Big(\{\widetilde{\pi}_{p}^{\varphi}\}_{p=1}^m, \widetilde{V}_{\varphi}, \{\widetilde{\mathcal{H}}^{\varphi}; \widetilde{\mathcal{E}}^{\varphi};\widetilde{\mathcal{D}}^{\varphi}\}\Big)$ with the following additional properties
  \begin{align}\label{Eq: Minimal condition}
  \widetilde{\mathcal{E^{\varphi}}} \subseteq \mathcal{E^{\varphi}} \text{~~and~~} \widetilde{\mathcal{H}_{\ell}^{\varphi}}
   =\Big[\big(\prod\limits_{p=1}^m \widetilde{\pi}_{p}^{\varphi}(\mathcal{A})\big)\widetilde{V}_{\varphi}\mathcal{H}_{\ell}\Big]\; \text{for every}\; \ell\in \Omega.
   \end{align}
   \end{prop}
   \begin{proof}First, we define a family of Hilbert spaces by $\widetilde{\mathcal{H}}_{\ell}^{\varphi}:= \big[\big(\prod\limits_{p=1}^m
   \pi^{\varphi}_{p}(\mathcal{A})\big)V_{\varphi}\mathcal{H}_{\ell}\big]$ for each $\ell\in \Omega$ and 
   $\widetilde{\mathcal{H}}^{\varphi}:= \big[\big(\prod\limits_{i=1}^m
   \pi^{\varphi}_{p}(\mathcal{A})\big)V_{\varphi}\mathcal{D}\big].$ Also take $ \widetilde{\mathcal{E}^{\varphi}}:= \{\widetilde{\mathcal{H}}^{\varphi}_{\ell}: \ell\in \Omega\}$ and $\widetilde{\mathcal{D}}^{\varphi}: =\bigcup\limits_{\ell\in \Omega}\widetilde{\mathcal{H}}^{\varphi}_{\ell}.$ Then by Lemma \ref{Lemma:  Technical}, it follows that $ \{\widetilde{\mathcal{H}}^{\varphi}; \widetilde{\mathcal{E}}^{\varphi};\widetilde{\mathcal{D}}^{\varphi}\}$ is a quantized domain. Since $\widetilde{\mathcal{H}}^{\varphi}_{\ell}\subseteq \mathcal{H}^{\varphi}_{\ell}$ for each $\ell\in\Omega,$ we have 
   $\widetilde{\mathcal{E}}^{\varphi}\subseteq \mathcal{E}^{\varphi}.$ For each $1\leq p\leq m$, we define a map $\widetilde{\pi}^{\varphi}_{p}: \mathcal{A}\to C^{*}_{\widetilde{\mathcal{E}}^{\varphi}}(\widetilde{\mathcal{D}}^{\varphi})$ by
   \begin{equation}\label{Equation: restriction}
       \widetilde{\pi}^{\varphi}_{p}(a)= {\pi}^{\varphi}_{p}(a)\big|_{\widetilde{\mathcal{D}}^{\varphi}}.
   \end{equation}
    Suppose $\sum\limits_{j=1}^n \prod\limits_{p=1}^{m}\pi_{p}^{\varphi}(a_{p_{j}})V_{\varphi}h_{j}\in \Span\{\prod\limits_{p=1}^m\pi^{\varphi}_{p}(\mathcal{A}))V_{\varphi}\mathcal{H}_{\ell}\},$ where $a_{p_{j}}\in \mathcal{A},$ and $h_{j}\in \mathcal{H}_{\ell}$ for  $1\leq j\leq n.$ Then by using the fact that the family $\{\pi^{\varphi}_{p}\}_{p=1}^{m}$ is commutative, we see that
    \begin{align*}
   \pi_{q}^{\varphi}(a)&\Big( \sum\limits_{j=1}^n \prod\limits_{p=1}^{m} \pi_{p}^{\varphi}(a_{p_{j}})V_{\varphi}h_{j}\Big)
  = \sum\limits_{j=1}^n\pi^{\varphi}_{1}(a_{1_{j}}) \cdots \pi_{q}^{\varphi}(a a_{q_{j}})\cdots \pi^{\varphi}_{m}(a_{m_{j}})V_{\varphi}h_{j} \in \widetilde{\mathcal{H}_{\ell}^{\varphi}}\end{align*}
  for each $1\leq q\leq m.$
  It follows that $$\pi_{p}^{\varphi}(a)(\Span\{ (\prod\limits_{p=1}^m \pi_{p}^{\varphi}(\mathcal{A}))V_{\varphi}\mathcal{H}^{\varphi}\})
  \subseteq \widetilde{\mathcal{H}}_{\ell}^{\varphi} \text{~for every~} 1\leq p\leq m.$$ 
  Since $\pi_{p}^{\varphi}(a)|_{\widetilde{\mathcal{H}}_{\ell}^{\varphi}}$ is continuous, wehave  $\pi_{p}^{\varphi}(a)(\widetilde{\mathcal{H}}_{\ell}^{\varphi})\subseteq \widetilde{\mathcal{H}}_{\ell}^{\varphi}.$ Now from the definition of $\widetilde{\pi}_{p}^{\varphi}$ we conclude that $\widetilde{\pi}_{p}^{\varphi}(a)(\widetilde{\mathcal{H}}_{\ell}^{\varphi})\subseteq \widetilde{\mathcal{H}}_{\ell}^{\varphi}$  for all $a\in \A$ and for each $\ell \in \Omega.$ It follows that $\widetilde{\pi}^{\varphi}_{p}(a)\in C^*_{\widetilde{\mathcal{E}}_{\varphi}}(\widetilde{\mathcal{D}}^{\varphi}).$ Moreover, $\{\widetilde{\pi}^{\varphi}_{p}\}$ is a commuting family and hence $\Big(\{\widetilde{\pi}_{p}^{\varphi}\}_{p=1}^m, \widetilde{V}_{\varphi}, \{\widetilde{\mathcal{H}}^{\varphi}; \widetilde{\mathcal{E}}^{\varphi};\widetilde{\mathcal{D}}^{\varphi}\}\Big)$ is a Stinespring's triple for $\varphi$ follows from Equation \eqref{Equation: restriction}. 
     \end{proof}
\begin{defn}\label{Definition: minimality}
    The Stinespring triple $\Big(\{\pi_{p}^{\varphi}\}_{p=1}^m, V_{\varphi}, \{\mathcal{H}^{\varphi}; \mathcal{E}^{\varphi};\mathcal{D}^{\varphi}\}\Big)$ of  $\varphi \in ICPCC^{s}_{loc}(\mathcal{A}^{k}, C^{\ast}_{\mathcal{E}}(\mathcal{D}))$ is called \emph{minimal}, if 
\begin{equation*}
    \mathcal{H}^{\varphi}_{\ell} = \big[\prod\limits_{p=1}^{m}\pi_{p}^{\varphi}(\mathcal{A})V_{\varphi}\mathcal{H}_{\ell}\big], \; \text{for each} \; \ell \in \Omega.
\end{equation*}
\end{defn}

Now, our aim is to show that the minimality condition in the Definition \ref{Definition: minimality} is appropriate in the sense that
there is a unitary equivalence relation between any two minimal Stinespring representations.  
\begin{thm}\label{Theorem:unitary equivalence of minimal Stinespring representation} Let $\Big(\{\widetilde{\pi}_{p}^{\varphi}\}_{p=1}^m, \widetilde{V}_{\varphi}, \{\widetilde{\mathcal{H}}^{\varphi}; \widetilde{\mathcal{E}}^{\varphi};\widetilde{\mathcal{D}}^{\varphi}\}\Big)$ and $\Big(\{\widehat{\pi}_{p}^{\varphi}\}_{p=1}^m, \widehat{V}_{\varphi}, \{\widehat{\mathcal{H}}^{\varphi}; \widehat{\mathcal{E}}^{\varphi};\widehat{\mathcal{D}}^{\varphi}\}\Big)$ be two mimimal Stinespring triples of the map $\varphi\in ICPCC^{s}_{loc}(\mathcal{A}^k, C^{*}_{\mathcal{E}}(\mathcal{D})).$ Then there exists a unitary operator $U_{\varphi}: \widetilde{\mathcal{H}}^{\varphi}\to \widehat{\mathcal{H}}^{\varphi}$ such that 
\begin{enumerate}
       \item $U_{\varphi}\widetilde{V}_{\varphi}=\widehat{V}_{\varphi}$ and 
    \item $U_{\varphi} \prod\limits_{p=1}^m \widetilde{\pi}_{p}^{\varphi}(a_{p})= \prod\limits_{p=1}^m\widehat{\pi}_{p}^{\varphi}(a_{p})U_{\varphi} $ for every $a_{p}\in \mathcal{A}$, $p\in \{1, \hdots, m\}.$
\end{enumerate}
\end{thm}
\begin{proof} By minimality, we have $ \big[\prod\limits_{p=1}^{m}\widetilde{\pi}_{p}^{\varphi}(\mathcal{A})\widetilde{V}_{\varphi}\mathcal{H}_{\ell}\big]=\widetilde{\mathcal{H}}^{\varphi}_{\ell}$ and $ \big[\prod\limits_{p=1}^{m}\widehat{\pi}_{p}^{\varphi}(\mathcal{A})\widehat{V}_{\varphi}\mathcal{H}_{\ell}\big] = \widehat{\mathcal{H}}^{\varphi}_{\ell}$
for each $\ell \in \Omega$. Then by Lemma \ref{Lemma:  Technical}, we get
$\big[\prod\limits_{p=1}^{m}\widetilde{\pi}_{p}^{\varphi}(\mathcal{A})\widetilde{V}_{\varphi}\mathcal{D}\big]=\widetilde{\mathcal{H}}^{\varphi}$
and $\big[\prod\limits_{p=1}^{m}\widehat{\pi}_{p}^{\varphi}(\mathcal{A})\widehat{V}_{\varphi}\mathcal{D}\big]=\widehat{\mathcal{H}}^{\varphi}.$
Let us define a map 
$U_{\varphi} \colon \Span\big\{\prod\limits_{p=1}^{m}\widetilde{\pi}_{p}^{\varphi}(\mathcal{A})\widetilde{V}_{\varphi}\mathcal{D}\big\} \to \Span\big\{\prod\limits_{p=1}^{m}\widehat{\pi}_{p}^{\varphi}(\mathcal{A})\widehat{V}_{\varphi}\mathcal{D}\big\}
$  by
\begin{equation*}
  U_{\varphi}\Big(\sum\limits_{j=1}^n \prod\limits_{p=1}^{m}\widetilde{\pi}_{p}^{\varphi}(a_{{p}_{j}}\big)
\widetilde{V}_{\varphi}h_{j}\Big)=  
\sum\limits_{j=1}^n\prod\limits_{p=1}^{m}\widehat{\pi}_{p}^{\varphi}(a_{{p}_{j}})\widehat{V}_{\varphi}h_{j} , 
\end{equation*}
 for all $a_{p_{j}}\in \mathcal{A}, \text{~and ~} h_{j}\in \mathcal{D}$. Firstly, we show that $U_{\varphi}$ is well defined. Suppose that $k = 2m-1$ (the result follows for $k = 2m$) and without loss of generality there is an $\ell \in\Omega$ such that $h_{j}\in \mathcal{H}_{\ell}$ for $j\in \{1,\hdots, n\}.$ It follows that 
\begin{align*}
  \Big\| U_{\varphi}\Big(\sum\limits_{j=1}^n\prod\limits_{p=1}^{m}\widetilde{\pi}_{p}^{\varphi}(a_{{p}_{j}})
\widetilde{V}_{\varphi}h_{j}\Big)\Big\|_{\widehat{\mathcal{H}}_{\ell}^{\varphi}}^2 & = \Big\|\sum\limits_{j=1}^n\prod\limits_{p=1}^{m}\widehat{\pi}_{p}^{\varphi}(a_{{p}_{j}})
\widehat{V}_{\varphi}h_{j}\Big\|^2_{\widehat{\mathcal{H}}^{\varphi}_{\ell}}\\
&=\sum\limits_{i,j=1}^n \Big\langle \prod\limits_{p=1}^{m}\widehat{\pi}_{p}^{\varphi}(a_{{p}_{j}})
\widehat{V}_{\varphi}h_{j}, \prod\limits_{p=1}^{m}\widehat{\pi}_{p}^{\varphi}(a_{{p}_{i}})
\widehat{V}_{\varphi}h_{i}\Big\rangle_{\widehat{\mathcal{H}}^{\varphi}_{\ell}}\\
&= \sum\limits_{i,j=1}^n
\big \langle (\widehat{V}_{\varphi})^* \prod\limits_{p=1}^{m}\widehat{\pi}_{p}^{\varphi}(a_{{p}_{i}}^*a_{{p}_{j}})
\widehat{V}_{\varphi}h_{j}, h_{i}\big \rangle_{{\mathcal{H}}_{\ell}}\\
&=\sum\limits_{i,j=1}^n \big \langle \varphi(a_{m_i}^*, \hdots,  a_{p_i}^*a^*, \hdots, a_{2_i}^*,a_{1_i}^*a_{1_j}, \hdots ,aa_{p_j}, \hdots, a_{m_j})h_{j}, h_{i}\big \rangle_{\mathcal{H}_{\ell}}\\
&= \sum\limits_{i,j=1}^n
\big \langle (\widetilde{V}_{\varphi})^* \prod\limits_{p=1}^{m}\widetilde{\pi}_{p}^{\varphi}(a_{{p}_{i}}^*a_{{p}_{j}})
\widetilde{V}_{\varphi}h_{j}, h_{i}\big \rangle_{{\mathcal{H}}_{\ell}}\\
&= \Big\| \sum\limits_{j=1}^n \prod\limits_{p=1}^{m}\widetilde{\pi}_{p}^{\varphi}(a_{{p}_{j}}\big)
\widetilde{V}_{\varphi}h_{j}\Big\|_{\widetilde{H}_{\ell}^{\varphi}}.
\end{align*}
It implies that the linear map
$U_{\varphi} \colon \Span\big\{\prod\limits_{p=1}^{m}\widetilde{\pi}_{p}^{\varphi}(\mathcal{A})\widetilde{V}_{\varphi}\mathcal{D}\big\} \to \Span\big\{\prod\limits_{p=1}^{m}\widehat{\pi}_{p}^{\varphi}(\mathcal{A})\widehat{V}_{\varphi}\mathcal{D}\big\}$ is a well defined and onto map.  Moreover $U_{\varphi}$ is an isometry.  Since $\Span\big\{\prod\limits_{p=1}^{m}\widehat{\pi}_{p}^{\varphi}(\mathcal{A})\widehat{V}_{\varphi}\mathcal{D}\big\}$ is dense in $\widehat{\mathcal{H}}^{\varphi}$ and 
$\Span\big\{\prod\limits_{p=1}^{m}\widetilde{\pi}_{p}^{\varphi}(\mathcal{A})\widetilde{V}_{\varphi}\mathcal{D}\big\}$ is dense in $\widetilde{\mathcal{H}}^{\varphi},$ thus $U_{\varphi}$ can be extended to a unitary from $\widetilde{\mathcal{H}}^{\varphi}$ to $\widehat{\mathcal{H}}^{\varphi}.$ We denote this unitary map  by $U_{\varphi}$ itself. By construction of $U_{\varphi}$, it is easy to see that $U_{\varphi}(\widetilde{\mathcal{H}}^{\varphi}_{\ell}) \subseteq \widehat{\mathcal{H}}^{\varphi}_{\ell}$  for all $\ell \in \Omega.$ 

Proof of (1): Let $h\in \mathcal{D}$.
 Then $U_{\varphi}\widetilde{V}_{\varphi}h=U_{\varphi}\big (\prod\limits_{p=1}^{m}\widetilde{\pi}_{p}^{\varphi}(1)\widetilde{V}_{\varphi}h \big)
 = \prod\limits_{p=1}^{m}\widehat{\pi}_{p}^{\varphi}(1)\widehat{V}_{\varphi}h =\widehat{V}_{\varphi}h.$ Since  $U_{\varphi}\widetilde{V}_{\varphi},\widehat{V}_{\varphi}$ are bounded linear operators which agree on $\mathcal{D}$, we conclude that $U_{\varphi}\widetilde{V}_{\varphi}=\widehat{V}_{\varphi}.$
 
  Proof of (2): $\Dom\big(U_{\varphi} \prod\limits_{p=1}^m \widetilde{\pi}_{p}^{\varphi}(a_{p})\big)=\{\xi\in \Dom(\prod\limits_{p=1}^m \widetilde{\pi}_{p}^{\varphi}(a_{p}))\colon \prod\limits_{p=1}^m \widetilde{\pi}_{p}^{\varphi}(a_{p})(\xi)\in \Dom(U_{\varphi})\}=
\widetilde{\mathcal{D}}^{\varphi},$ and
$\Dom(\prod\limits_{p=1}^m\widehat{\pi}_{p}^{\varphi}(a_{p})U_{\varphi})=\{\xi \in \widetilde{\mathcal{H}}^{\varphi}\colon U_{\varphi}(\xi)\in \Dom(\prod\limits_{p=1}^m\widehat{\pi}_{p}^{\varphi}(a_{p}))\}=\{\xi \in \widetilde{\mathcal{H}}^{\varphi}\colon U_{\varphi}(\xi)\in  \widehat{\mathcal{D}}^\varphi\}=\widetilde{\mathcal{D}}^\varphi.$ Suppose $\xi$ is of the form $\xi=\sum\limits_{j=1}^n\prod\limits_{p=1}^{m}\widetilde{\pi}_{p}^{\varphi}(a_{{p}_{j}})
\widetilde{V}_{\varphi}h_{j},$ where $a_{p_{j}}\in \mathcal{A}, h_{j}\in \mathcal{H}_{\ell}.$ Then for all $b_{p}\in \mathcal{A},$
we have
\begin{align*}
   U_{\varphi} \prod\limits_{p=1}^m 
  \widetilde{\pi}_{p}^{\varphi}(b_{p})\xi=U_{\varphi}\Big(\sum\limits_{j=1}^n\prod\limits_{p=1}^{m}\widetilde{\pi}_{p}^{\varphi}(b_{p}a_{{p}
  _{j}}) \widetilde{V}_{\varphi}h_{j}\Big)
  &=\sum\limits_{j=1}^n\prod\limits_{p=1}^{m}\widehat{\pi}_{p}^{\varphi}(b_{p}a_{{p}_{j}}) \widehat{V}_{\varphi}h_{j}\\
  &=\prod\limits_{p=1}^{m}\widehat{\pi}_{p}^{\varphi}(b_{p})\Big(\sum\limits_{j=1}^n\prod\limits_{p=1}^{m}\widehat{\pi}_{p}^{\varphi}(a_{{p}_{j}})   \widehat{V}_{\varphi}h_{j}\Big)\\
  &=\prod\limits_{p=1}^{m}\widehat{\pi}_{p}^{\varphi}(b_{p})U_{\varphi} \xi.
\end{align*}
This implies that
$U_{\varphi} \prod\limits_{p=1}^m \widetilde{\pi}_{p}^{\varphi}(b_{p})|_{\widetilde{\mathcal{H}}^{\varphi}_{\ell}}= \prod\limits_{p=1}^m\widehat{\pi}_{p}^{\varphi}(b_{p})U_{\varphi}|_{\widetilde{\mathcal{H}}^{\varphi}_{\ell}}$ for all $\ell\in \Omega.$ This completes the proof. 
\end{proof}
\section{Radon-Nikod\'{y}m theorem}
In this section, we prove Radon-Nikod\'{y}m type theorem  for  local completely positive invariant $k$-linear maps. In $C^{\ast}$-algebra theory, a Radon-Nikod\'{y}m derivative is obtained between two CP-maps that are dominated by one with the other map. So, in this setup, it is necessary to define such relation $``\leq"$ on $ ICPCC^{s}_{loc}(\mathcal{A}^{k}, C^{\ast}_{\mathcal{E}}(\mathcal{D})).$ The results proved in this section generalizes our earlier work with  B. V. Rajarama Bhat \cite{BGS} when $k=1$.
\begin{defn} Let $\varphi, \psi \in ICPCC^{s}_{loc}(\mathcal{A}^{k}, C^{\ast}_{\mathcal{E}}(\mathcal{D})).$ If $\varphi-
\psi\in  ICPCC^{s}_{loc}(\mathcal{A}^{k}, C^{\ast}_{\mathcal{E}}(\mathcal{D})),$ then $\psi$ is said to be \emph{dominated} by $\varphi.$ In this case, we write $\psi\leq \varphi.$ 
\end{defn}
Let $\varphi \in ICPCC^{s}_{loc}(\mathcal{A}^{k}, C^{\ast}_{\mathcal{E}}(\mathcal{D}))$ and let $\Big(\{\pi_{p}^{\varphi}\}_{p=1}^m, V_{\varphi}, \{\mathcal{H}^{\varphi}; \mathcal{E}^{\varphi};\mathcal{D}^{\varphi}\}\Big)$ be the minimal Stinespring triple of $\varphi.$ The commutant of the family $\{\pi_{p}^{\varphi}\}_{p=1}^m$ is defined by the subspace of  $\mathcal{B}(\mathcal{H}^{\varphi})$ as,
\begin{align} 
\bigcap_{p=1}^m \pi_{p}^{\varphi}(\mathcal{A})'= \{T\in\mathcal{B}(\mathcal{H}^{\varphi}): T \pi^{\varphi}_{p}(a)\subseteq \pi_{p}^{\varphi}(a)T \text{~ for all~} a \in \mathcal{A}, 1\leq p\leq m   \}.
\end{align}
\begin{lemma}\label{Lemma: Existenceofboundedoperator}
Let $\varphi_{1}, \varphi_{2} \in ICPCC^{s}_{loc}(\mathcal{A}^{k}, C^{\ast}_{\mathcal{E}}(\mathcal{D}))$ and $m = \big[\frac{k+1}{2}\big]$. Let $\Big( \{\pi_{p}^{\varphi_{i}}\}_{p=1}^{m}, V_{\varphi_{i}}, \{\mathcal{H}^{\varphi_{i}}; \mathcal{E}^{\varphi_{i}}; \mathcal{D}^{\varphi_{i}}\}\Big)$ be the minimal Stinespring's triple for $\varphi_{i}$ with $i =1,2$. If $\varphi_{1} \leq \varphi_{2}$, then there exists a unique contraction $T \in \mathcal{B}(\mathcal{H}^{\varphi_{2}}, \mathcal{H}^{\varphi_{1}})$ such that 
\begin{enumerate}
    \item[(i)] $TV_{\varphi_{2}} = V_{\varphi_{1}} $.
    \item[(ii)]$T\Big(\prod\limits_{p=1}^{m}\pi^{\varphi_{2}}_{p}(a_{p})\Big) \subseteq \Big(\prod\limits_{p=1}^{m}\pi^{\varphi_{1}}_{p}(a_{p})\Big)T$, for all $a_{p}\in \mathcal{A}$ and $p = 1,2,\hdots, m$. 
\end{enumerate}
\end{lemma}
\begin{proof} It is given that  $(\varphi_{2}-\varphi_{1}) \in ICPCC_{loc}(\mathcal{A}^{k}, C^{\ast}_{\mathcal{E}}(\mathcal{D}))$. Let us take  ${D}_{p}$ (for $2 \leq p \leq m$) and  $R_{1}$  as in  Equation \eqref{Eq: Essential}.  Suppose that $k = 2m-1$, then
 $(D_{m}^{\ast}, \hdots, D_{2}^{\ast}, R_{1}^{\ast}R_{1}, D_{2}, \hdots, D_{m})\in M_{n}(\mathcal{A})^{2m-1}$ is $\alpha$-symmetric (see Definition \ref{Definition: alpha-symmetric}) and $R_{1}^{\ast}R_{1} \geq_{\alpha} 0$, for every $\alpha \in \Lambda$. Moreover, for $a_{p_{j}}\in \mathcal{A},\; h_{j} \in \mathcal{H}_{\ell} $ with $j\in \{1,\hdots,n\}$ and $p\in \{1, \hdots, m\},\; n \in \mathbb{N}$ we see that 
 \begin{align}\label{Equation: pi2dominatespi1}
   \notag  \sum\limits_{t,j = 1}^{n}\Big\langle \varphi_{1}(a_{m_t}^{\ast}, \hdots,a_{2_t}^{\ast}, a_{1_t}^{\ast}a_{1_{j}}, a_{2_{j}},& \hdots, a_{m_{j}})h_{j},\; h_{t} \Big\rangle_{\mathcal{H}_{\ell}}\\
   \notag  &= \Big\langle (\varphi_{1})_{n} (D_{m}^{\ast},\hdots, D_{2}^{\ast}, R_{1}^{\ast}R_{1}, D_{2}, \hdots, D_{m}) \begin{bmatrix}
         h_{1}\\ \vdots \\ h_{n}
         \end{bmatrix}, \; \begin{bmatrix}
         h_{1}\\ \vdots \\ h_{n}
         \end{bmatrix}
         \Big\rangle_{\mathcal{H}_{\ell}^{n}}\\
       \notag  &\leq \Big\langle (\varphi_{2})_{n} (D_{m}^{\ast},\hdots, D_{2}^{\ast}, R_{1}^{\ast}R_{1}, D_{2}, \hdots, D_{m}) \begin{bmatrix}
         h_{1}\\ \vdots \\ h_{n}
         \end{bmatrix}, \; \begin{bmatrix}
         h_{1}\\ \vdots \\ h_{n}
         \end{bmatrix}
         \Big\rangle_{\mathcal{H}_{\ell}^{n}}\\
         &= \sum\limits_{t,j = 1}^{n}\Big\langle \varphi_{2}(a_{m_t}^{\ast}, \hdots,a_{2_t}^{\ast}, a_{1_t}^{\ast}a_{1_{j}}, a_{2_{j}}, \hdots, a_{m_{j}})h_{j},\; h_{t} \Big\rangle_{\mathcal{H}_{\ell}}.
 \end{align}
The above observation is helpful in defining the desired contraction from $\mathcal{H}^{\varphi_{2}}$ to $\mathcal{H}^{\varphi_{1}}$. Firstly, by the minimaility condition from Definition \ref{Definition: minimality}, we have 
\begin{equation*}
    \mathcal{H}_{\ell}^{\varphi_{i}} = \Big[ \prod\limits_{p=1}^{m}\pi_{p}^{\varphi_{i}}(\mathcal{A})V_{\varphi_{i}}\mathcal{H}_{\ell}\Big],\; \text{for every}\; \ell \in \Omega, \; i = 1,2.
\end{equation*}
Thus by Lemma \ref{Lemma: Technical}, it implies that $   \mathcal{H}^{\varphi_{i}} = \Big[ \Big( \prod\limits_{p=1}^{m} \pi_{p}^{\varphi_{i}}(\mathcal{A})\Big)V_{\varphi_{i}}\mathcal{D}\Big] , \; \text{for}\; i = 1,2.$ Let us define a map $T$ on the dense subset of $\mathcal{H}^{\varphi_{2}}$ as,
\begin{equation}\label{Equation: definitionofT}
    T\Bigg( \sum\limits_{j=1}^{n}\big(\prod\limits_{p=1}^{m} \pi_{p}^{\varphi_{2}}({a_{p_{j}}})V_{\varphi_{2}}h_{j}\big)\Bigg) = \sum\limits_{j=1}^{n}\big(\prod\limits_{p=1}^{m} \pi_{p}^{\varphi_{1}}({a_{p_{j}}})V_{\varphi_{1}}h_{j}\big),
\end{equation}
for all $a_{p_{j}} \in \mathcal{A}$, $h_{j} \in \mathcal{D}$, $j = \{1, 2, \hdots n\}$ and $p = \{1, 2, \hdots, m\}$. Now we show that $T$ is a well defined bounded map. Since $\mathcal{D}$ is the union space, there is an $\ell \in \Omega$ with $h_{j}\in \mathcal{H}_{\ell}$, for $j \in \{1, 2, \hdots, n\}$. Then we have 
\begin{align*}
    \Bigg\|\sum\limits_{j=1}^{n}\big(\prod\limits_{p=1}^{m} \pi_{p}^{\varphi_{1}}({a_{p_{j}}})V_{\varphi_{1}}h_{j}\big) \Bigg\|_{\mathcal{H}_{\ell}^{\varphi_{1}}}^{2}& =\Bigg\langle \sum\limits_{j=1}^{n}\big(\prod\limits_{p=1}^{m} \pi_{p}^{\varphi_{1}}({a_{p_{j}}})V_{\varphi_{1}}h_{j}\big),\; \sum\limits_{j=1}^{n}\big(\prod\limits_{p=1}^{m} \pi_{p}^{\varphi_{1}}({a_{p_{j}}})V_{\varphi_{1}}h_{j}\big) \Bigg\rangle_{\mathcal{H}_{\ell}^{\varphi_{1}}}\\
    &= \sum\limits_{t, j= 1}^{n}\Big\langle V_{\varphi_{1}}^{\ast}\pi_{1}^{\varphi_{1}}(a_{1_{t}}^{\ast})\cdots\pi_{m}^{\varphi_{1}}(a_{m_{t}}^{\ast})\pi_{1}^{\varphi_{1}}(a_{m_{j}})\cdots \pi_{1}^{\varphi_{1}}(a_{1_{j}})V_{\varphi_{1}}h_{j},\; h_{t} \Big\rangle_{\mathcal{H}_{\ell}}\\
    &= \sum\limits_{t, j= 1}^{n}\Big\langle V_{\varphi_{1}}^{\ast}\pi_{1}^{\varphi_{1}}(a_{1_{t}}^{\ast}a_{1_{j}})\cdots\pi_{2}^{\varphi_{1}}(a_{2_{t}}^{\ast}a_{2_{j}})\cdots \pi_{m}^{\varphi_{1}}(a_{m_{t}}^{\ast}a_{m_{j}})V_{\varphi_{1}}h_{j},\; h_{t} \Big\rangle_{\mathcal{H}_{\ell}}\\
    &=\sum\limits_{t, j= 1}^{n}\Big\langle \varphi_{1}(a_{m_{t}}^{\ast}, \hdots,a_{2_{t}}^{\ast}, a_{1_{t}}^{\ast}a_{1_{j}}, a_{2_{j}}, \hdots, a_{m_{j}})h_{j},\; h_{t} \Big\rangle_{\mathcal{H}_{\ell}}.
\end{align*}
Further by Equation \eqref{Equation: pi2dominatespi1}, we get 
\begin{equation*}
    \Big\|  T\Bigg( \sum\limits_{j=1}^{n}\big(\prod\limits_{p=1}^{m} \pi_{p}^{\varphi_{2}}({a_{p_{j}}})V_{\varphi_{2}}h_{j}\big)\Bigg)\Big\| \leq  \Big\|\sum\limits_{j=1}^{n}\big(\prod\limits_{p=1}^{m} \pi_{p}^{\varphi_{1}}({a_{p_{j}}})V_{\varphi_{1}}h_{j}\big)\Big\|.
\end{equation*}
This shows that $T$ is a well defined contraction mapping from a dense subsapce of $\mathcal{H}^{\varphi_{2}}$ onto the dense subspace of $\mathcal{H}^{\varphi_{1}}$. Thus there is a unique contraction defined on $\mathcal{H}^{\varphi_{2}}$, again denote it by $T \in \mathcal{B}(\mathcal{H}^{\varphi_{2}}, \mathcal{H}^{\varphi_{1}})$. It is clear from the definition $T$ in Equation \eqref{Equation: definitionofT} that $T(\mathcal{H}_{\ell}^{\varphi_{2}}) \subseteq \mathcal{H}_{\ell}^{\varphi_{1}}$ and $T\big|_{\mathcal{H}_{\ell}^{\varphi_{2}}} \in \mathcal{B}(\mathcal{H}_{\ell}^{\varphi_{2}}, \mathcal{H}_{\ell}^{\varphi_{1}})$, for every $\ell \in \Omega.$

 \noindent Proof of (i): Let $g \in \mathcal{D}$. Then
\begin{equation*}
    TV_{\varphi_{2}}g = T \Big( \prod\limits_{p=1}^{m}\pi_{p}^{\varphi_{2}}(1)V_{\varphi_{2}}g\Big) = \prod\limits_{p=1}^{m}\pi_{p}^{\varphi_{1}}(1)V_{\varphi_{1}}g = V_{\varphi_{1}}g.
\end{equation*}
The result follows from the fact that $\mathcal{D}$ is dense in $\mathcal{H}$.

\noindent Proof of (ii): 
Let $a_{p}, b_{p} \in \mathcal{A}$ with $p \in \{1, \hdots, m\}$ and $h\in \mathcal{D}$. Then by using the fact that $\{\pi_{p}^{\varphi_{i}}\}_{p=1}^{m}$ is a commuting family of $\ast$-homomorphisms from $\mathcal{A}$ to $C^{\ast}_{\mathcal{E}^{\varphi_{i}}}(\mathcal{D}^{\varphi_{i}})$ for $i = 1,2$ we have
\begin{align*}
    T\Big(\prod\limits_{p=1}^{m}\pi^{\varphi_{2}}_{p}(a_{p})\Big) \Big( \prod\limits_{p=1}^{m} \pi_{p}^{\varphi_{2}}(b_{p})V_{\varphi_{2}}h\Big) &= T\Big(\prod\limits_{p=1}^{m}\pi^{\varphi_{2}}_{p}(a_{p}b_{p})V_{\varphi_{2}}h\Big)  \\
    &= \prod\limits_{p=1}^{m}\pi^{\varphi_{1}}_{p}(a_{p}b_{p})V_{\varphi_{1}}h\\
    &= \prod\limits_{p=1}^{m}\pi^{\varphi_{1}}_{p}(a_{p})\Big( \prod\limits_{p=1}^{m} \pi_{p}^{\varphi_{2}}(b_{p})V_{\varphi_{1}}h\Big)\\
    &= \prod\limits_{p=1}^{m}\pi^{\varphi_{1}}_{p}(a_{p})T\Big( \prod\limits_{p=1}^{m} \pi_{p}^{\varphi_{2}}(b_{p})V_{\varphi_{2}}h\Big).
\end{align*}
It follows that $T\Big(\prod\limits_{p=1}^{m}\pi^{\varphi_{2}}_{p}(a_{p})\Big) \subseteq \Big(\prod\limits_{p=1}^{m}\pi^{\varphi_{1}}_{p}(a_{p})\Big)T$, for all $a_{p}\in \mathcal{A}$ and $p \in \{1, \hdots, m\}$. Similarly, if $k = 2m$, then $(D_{m}^{\ast}, \hdots, D_{2}^{\ast}, R_{1}^{\ast}, R_{1}, D_{2}, \hdots, D_{m}) \in \mathcal{A}^{2m}$ is $\alpha$-symmetric, for every $\alpha \in \Lambda$ and hence the result follows same as in the odd case. 
\end{proof}
\begin{lemma} \label{Lemma: introductionofphiT} Let $\varphi\in ICPCC^{s}_{loc}(\mathcal{A}^{k}, C^{\ast}_{\mathcal{E}}(\mathcal{D}))$ and $m = \big[\frac{k+1}{2}\big]$. For every $T \in \bigcap\limits_{p=1}^{m}\pi^{\varphi}_{p}(\mathcal{A})^{\prime}\cap C^{\ast}_{\mathcal{E}^{\varphi}}(\mathcal{D}^{\varphi})$, define the $k$-linear map $\varphi_{T}\colon \mathcal{A}^{k} \to C_{\mathcal{E}}(\mathcal{D})$ by 
\begin{align}\label{Eq: varphiT}
\varphi_{T}(a_{1}, \hdots, a_{k}) = \left\{ \begin{array}{cc}
     V_{\varphi}^{\ast}T\pi^{\varphi}_{1}(a_{m}) \pi^{\varphi}_{2}(a_{m-1}a_{m+1})\cdots \pi_{m}^{\varphi}(a_{1}a_{2m-1})V_{\varphi},&\text{if} \; k = 2m-1 \\
   &\\
     V_{\varphi}^{\ast}T\pi^{\varphi}_{1}(a_{m}a_{m+1}) \pi^{\varphi}_{2}(a_{m-1}a_{m+2})\cdots \pi_{m}^{\varphi}(a_{1}a_{2m})V_{\varphi},&\text{if}\; k = 2m,
\end{array} \right.
\end{align}
for all $a_{1}, a_{2}, \hdots, a_{k} \in \mathcal{A}$. Then the correspondence $T \mapsto \varphi_{T}$ is an injective linear map .

\end{lemma}
\begin{proof}
Let $T_{1}, T_{2} \in \bigcap\limits_{p=1}^{m}\pi^{\varphi}_{p}(\mathcal{A})^{\prime}\cap C^{\ast}_{\mathcal{E}^{\varphi}}(\mathcal{D}^{\varphi})$ and $\lambda \in \mathbb{C}$. Then 
\begin{align*}
    \varphi_{(T_{1}+ \lambda T_{2})}&(a_{1}, \hdots, a_{k})\\
    &= \left\{ \begin{array}{cc}
     V_{\varphi}^{\ast}(T_{1}+\lambda T_{2})\pi^{\varphi}_{1}(a_{m}) \pi^{\varphi}_{2}(a_{m-1}a_{m+1})\cdots \pi_{m}^{\varphi}(a_{1}a_{2m-1})V_{\varphi},&\text{if} \; k = 2m-1 \\
   &\\
     V_{\varphi}^{\ast}(T_{1}+\lambda T_{2})\pi^{\varphi}_{1}(a_{m}a_{m+1}) \pi^{\varphi}_{2}(a_{m-1}a_{m+2})\cdots \pi_{m}^{\varphi}(a_{1}a_{2m})V_{\varphi},&\text{if}\; k = 2m
\end{array} \right.\\
   &= \varphi_{T_{1}}(a_{1},\hdots,a_{k}) + \lambda\;  \varphi_{T_{2}}(a_{1},\hdots,a_{k}).
\end{align*}
Suppose that $\varphi_{T} = 0$ for some $T \in \bigcap\limits_{p=1}^{m}\pi^{\varphi}_{p}(\mathcal{A})^{\prime}\cap C^{\ast}_{\mathcal{E}^{\varphi}}(\mathcal{D}^{\varphi})$. We claim that $T=0$. Let us take $a_{p_{j}}, b_{p_{i}}\in \mathcal{A}$ and $ h_{j}, g_{i} \in \mathcal{D}$, where $j\in \{1, \hdots, n\}, i \in \{1, \hdots, k\}$ and $p \in \{1, \hdots, m\}$. Then there exists an $\ell \in \Omega$ such that each $g_{j}, h_{i} \in \mathcal{H}_{\ell}$ for all $j\in \{1, \hdots, n\}, i \in \{1, \hdots, k\}.$ Then  we compute that 
\begin{align*}
    \Big\langle T\Big( \sum\limits_{j=1}^{n}&\prod\limits_{p=1}^{m}\pi^{\varphi}_{p}(a_{p_{j}})V_{\varphi}h_{j}\Big), \; \Big( \sum\limits_{i=1}^{t}\prod\limits_{p=1}^{m}\pi^{\varphi}_{p}(b_{p_{i}})V_{\varphi}g_{i}\Big)\Big\rangle_{\mathcal{H}_{\ell}^{\varphi}}\\
    &=\sum\limits_{j=1}^{n}\sum\limits_{i=1}^{t}\Big\langle V_{\varphi}^{\ast}T \prod\limits_{p=1}^{m}\pi^{\varphi}_{p}(b_{p_{i}}^{\ast}a_{p_{j}})V_{\varphi}h_{j}, \; g_{i} \Big\rangle_{\mathcal{H}_{\ell}} \\
    &=\sum\limits_{j=1}^{n}\sum\limits_{i=1}^{t}\Big\langle V_{\varphi}^{\ast}T\pi^{\varphi}_{1}(b_{1_{i}}^{\ast}a_{1_{j}}) \pi^{\varphi}_{2}(b_{2_{i}}^{\ast}a_{2_{j}})\cdots \pi^{\varphi}_{m}(b_{m_{i}}^{\ast}a_{m_{j}})V_{\varphi}h_{j}, \; g_{i}\Big\rangle_{\mathcal{H}_{\ell}}\\
    &= \left\{ \begin{array}{cc}
         \sum\limits_{j=1}^{n}\sum\limits_{i=1}^{t} \Big\langle \varphi_{T}(b_{m_{i}}^{\ast}, \hdots, b_{2_{i}}^{\ast}, b_{1_{i}}^{\ast}a_{1_{j}}, a_{2_{j}}, \hdots, a_{m_{j}})h_{j},\; g_{i}\Big\rangle & \; \text{if}\; k = 2m-1 \\
         & \\
         \sum\limits_{j=1}^{n}\sum\limits_{i=1}^{t} \Big\langle \varphi_{T}(b_{m_{i}}^{\ast}, \hdots, b_{2_{i}}^{\ast}, b_{1_{i}}^{\ast}, a_{1_{j}}, a_{2_{j}}, \hdots, a_{m_{j}})h_{j},\; g_{i}\Big\rangle & \; \text{if}\; k = 2m
    \end{array}\right.\\
    &= 0.
\end{align*}
Since $\mathcal{H}^{\varphi} = \big[ \prod\limits_{p=1}^{m}\pi^{\varphi}_{p}(\mathcal{A})V_{\varphi}\mathcal{D}\big]$, we conclude that $T=0$.
\end{proof}
\begin{prop}\label{Proposition: phiT}
 Let $\Big( \{\pi_{p}^{\varphi}\}_{p=1}^{m}, V_{\varphi}, \{\mathcal{H}^{\varphi}; \mathcal{E}^{\varphi}; \mathcal{D}^{\varphi}\}\Big)$ be the minimal Stinespring's triple for  $\varphi \in ICPCC_{loc}^{s}(\mathcal{A}^{k}, C^{\ast}_{\mathcal{E}^{\varphi}}(\mathcal{D}^{\varphi}))$ and  $T \in \bigcap\limits_{p=1}^{m}\pi^{\varphi}_{p}(\mathcal{A})^{\prime}\cap C^{\ast}_{\mathcal{E}^{\varphi}}(\mathcal{D}^{\varphi})$. Let  $\varphi_{T} \colon \mathcal{A}^{k} \to C^{\ast}_{\mathcal{E}}(\mathcal{D}) $  be a map defined in Equation \eqref{Eq: varphiT}. 
 If $0 \leq T \leq I$, then the following assertions hold true.
\begin{enumerate}
    \item[(i)] $\varphi_{T} \in ICPCC_{loc}^{s}(\mathcal{A}^{k}, C^{\ast}_{\mathcal{E}}(\mathcal{D}))$.
     \item[(ii)] $\varphi_{T} \leq \varphi$.
\end{enumerate}
\end{prop}
\begin{proof} 

\noindent Proof of (i): For every $a_{1}, \hdots a_{k}\in \mathcal{A}$, we see that 
\begin{align*}
  &  \varphi^{\ast}_{T}(a_{1}, \hdots, a_{k}) \\
    &= \varphi_{T}(a_{k}^{\ast}, \hdots, a_{1}^{\ast})^{\ast}\\
    &= \left\{ \begin{array}{cc}
         V_{\varphi}^{\ast} \pi_{m}^{\varphi}(a_{2m-1}^{\ast}a_{1}^{\ast})^{\ast} \cdots \pi_{1}^{\varphi}(a_{m}^{\ast})^{\ast}TV_{\varphi}&\text{if}\; k=2m-1  \\
         & \\
         V_{\varphi}^{\ast} \pi_{m}^{\varphi}(a_{2m}^{\ast}a_{1}^{\ast})^{\ast} \cdots \pi_{1}^{\varphi}(a_{m+1}^{\ast}a_{m}^{\ast})^{\ast}TV_{\varphi}&\text{if}\; k=2m
    \end{array}
    \right.\\
    &\\
    &= \left\{ \begin{array}{cc}
     V_{\varphi}^{\ast}T\pi^{\varphi}_{1}(a_{m}) \cdots \pi_{m}^{\varphi}(a_{1}a_{2m-1})V_{\varphi},&\text{if} \; k = 2m-1 \\
   &\\
     V_{\varphi}^{\ast}T\pi^{\varphi}_{1}(a_{m}a_{m+1}) \cdots \pi_{m}^{\varphi}(a_{1}a_{2m})V_{\varphi},&\text{if}\; k = 2m
\end{array} \right.\\
&=  \varphi_{T}(a_{1}, \hdots, a_{k}).
\end{align*}
Since $\{\pi_{p}^{\varphi}\}_{p=1}^{m}$ is a commuting family of $\ast$-homomorphisms and $T \in \bigcap\limits_{p=1}^{m}\pi_{p}(\mathcal{A})^{\prime}\cap C^{\ast}_{\mathcal{E}}(\mathcal{D})$ is positive. This implies that $\varphi_{T}^{\ast} = \varphi_{T}$ and hence $\varphi_{T}$ is symmetric. 

Suppose $k = 2m-1$ and let $a_{1}, a_{2}, \hdots, a_{2m-1}, c_{1}, \hdots, c_{m-1} \in \mathcal{A}$. Then 
\begin{align*}
    \varphi_{T}(a_{1}c_{1}, \hdots, a_{m-1}c_{m-1}, &a_{m},a_{m+1},\hdots a_{2m-1})\\
    &= V_{\varphi}^{\ast}T\pi_{1}^{\varphi}(a_{m})\pi^{\varphi}_{2}(a_{m-1}c_{m-1}a_{m+1})\cdots \pi^{\varphi}_{m}(a_{1}c_{1}a_{2m-1})V_{\varphi}\\
    &=\varphi_{T}(a_{1}, \hdots, a_{m-1}, a_{m},c_{m-1}a_{m+1},\hdots c_{1}a_{2m-1}).
\end{align*}
Similarly, if $k = 2m$, we have     
\begin{equation*}
        \varphi_{T}(a_{1}c_{1}, \hdots,  a_{m-1}c_{m-1}, a_{m}c_{m},a_{m+1},\hdots a_{2m})=\varphi_{T}(a_{1}, \hdots, a_{m-1}, a_{m},c_{m}a_{m+1},\hdots c_{1}a_{2m}).
\end{equation*} 
Therefore, $\varphi_{T}$ is an invariant $k$-linear map. Now we show that $\varphi_{T}$ is local completely contractive. Suppose that $k=2m-1$. Let $n \in \mathbb{N}$ and $A_{t} = \big[a_{t_{ij}}\big]_{i,j=1}^{n} \in M_{n}(\A)^{2m-1}$ for $t = 1, \hdots, (2m-1)$. We know that the amplification map $(\varphi_{T})_{n} \colon M_{n}(\mathcal{A})^{k} \to C^{\ast}_{\mathcal{E}^{n}}(\mathcal{D}^{n})$ is defined by 
    \begin{align*}
        (\varphi_{T})_{n}&(A_{1}, A_{2}, \hdots, A_{2m-1})\\
        &= \Bigg[
         \sum\limits_{r_{1}, r_{2}, \hdots, r_{2(m-1)} = 1}^{n} \varphi_{T}\big(a_{1_{ir_{1}}}, a_{2_{r_{1}r_{2}}}, \hdots, a_{2m-1_{r_{2(m-1)}j}}\big)\Bigg]_{i,j=1}^{n}. 
    \end{align*}
    For every $h_{1}, h_{2}, \hdots, h_{n} \in \mathcal{D}$, due to the fact that $\mathcal{D}$ is the union space of upper filtered family $\mathcal{E}$, there exists an $\ell \in \Omega$ such that  $h_{1}, h_{2}, \hdots, h_{n} \in \mathcal{H}_{\ell}$. Then we have 
    \begin{align*}
        &\Bigg\|(\varphi_{T})_{n}(A_{1}, A_{2}, \hdots, A_{2m-1}) \begin{bmatrix}h_{1} \\ \vdots\\h_{n} \end{bmatrix} \Bigg\|_{\mathcal{H}_{\ell}^{n}}^{2} \\
        &= \sum\limits_{i=1}^{n}\Big\|\sum\limits_{j=1}^{n}\sum\limits_{r_{1}, r_{2}, \hdots, r_{2(m-1)} = 1}^{n} \varphi_{T}\big(a_{1_{ir_{1}}}, a_{2_{r_{1}r_{2}}}, \hdots, a_{2m-1_{r_{2(m-1)}j}}\big)h_{j}\Big\|_{\mathcal{H}_{\ell}}^{2}\\
        &= \sum\limits_{i=1}^{n}\Big\|\sum\limits_{j=1}^{n}\sum\limits_{r_{1}, r_{2}, \hdots, r_{2(m-1)} = 1}^{n} V_{\varphi}^{\ast}T \pi_{1}^{\varphi}(a_{m_{r_{m-1}}r_{m}}) \cdots \pi_{m}^{\varphi}(a_{1_{ir_{1}}}a_{{2m-1}_{r_{2(m-1)}j}})V_{\varphi}h_{j}\Big\|_{\mathcal{H}_{\ell}}^{2}\\
        &= \sum\limits_{i=1}^{n}\Big\|V_{\varphi}^{\ast}T\Big(\sum\limits_{j=1}^{n}\sum\limits_{r_{1}, r_{2}, \hdots, r_{2(m-1)} = 1}^{n}  \pi_{1}^{\varphi}(a_{m_{r_{m-1}}r_{m}}) \cdots \pi_{m}^{\varphi}(a_{1_{ir_{1}}}a_{{2m-1}_{r_{2(m-1)}j}})V_{\varphi}h_{j}\Big)\Big\|_{\mathcal{H}_{\ell}}^{2}\\
        \end{align*}
        \begin{align*}
        &\leq \sum\limits_{i=1}^{n}\Big\|\sum\limits_{j=1}^{n}\sum\limits_{r_{1}, r_{2}, \hdots, r_{2(m-1)} = 1}^{n}  \pi_{1}^{\varphi}(a_{m_{r_{m-1}}r_{m}}) \cdots \pi_{m}^{\varphi}\big(a_{1_{ir_{1}}}a_{{2m-1}_{r_{2(m-1)}j}}\big)V_{\varphi}h_{j}\Big\|_{\mathcal{H}^{\varphi}_{\ell}}^{2}\\
                &= \sum\limits_{i=1}^{n}\Big\|\sum\limits_{j=1}^{n}\sum\limits_{r_{1}, r_{2}, \hdots, r_{2(m-1)} = 1}^{n} \pi_{m}^{\varphi}(a_{1_{ir_{1}}}) \cdots \pi_{1}^{\varphi}(a_{m_{r_{m-1}}r_{m}})\cdots   \pi_{m}^{\varphi}\big(a_{{2m-1}_{r_{2(m-1)}j}}\big)V_{\varphi}h_{j}\Big\|_{\mathcal{H}^{\varphi}_{\ell}}^{2}\\
        &= \Bigg\| \big[\pi_{m}^{\varphi}(a_{1_{ij}})\big]_{i,j=1}^{n}\cdots \big[\pi_{1}^{\varphi}(a_{m_{ij}})\big]_{i,j=1}^{n}\cdots \big[\pi_{m}^{\varphi}(a_{2m-1_{ij}})\big]_{i,j=1}^{n} \begin{bmatrix}V_{\varphi}h_{1} \\ \vdots\\V_{\varphi}h_{n} \end{bmatrix}\Bigg\|_{{\mathcal{H}^{\varphi}_{\ell}}^{\oplus n}}^{2}\\
        & \leq \Big\|\big[\pi_{m}^{\varphi}(a_{1_{ij}})\big]_{i,j=1}^{n}\Big\|_{\ell}^{2} \cdots \Big\|\big[\pi_{1}^{\varphi}(a_{m_{ij}})\big]_{i,j=1}^{n}\Big\|^{2}_{\ell} \cdots \Big\|\big[\pi_{m}^{\varphi}(a_{2m-1_{ij}})\big]_{i,j=1}^{n} \Big\|_{\ell}^{2}\Bigg\|\begin{bmatrix}h_{1} \\ \vdots\\h_{n} \end{bmatrix} \Bigg\|_{{\mathcal{H}_{\ell}}^{\oplus n}}^{2}
    \end{align*}
    Since $\pi_{p}^{\varphi}$ is a local completely contractive map for $p = 1,\hdots, m$, there is an $\alpha_{p}\in \Lambda$ such that 
    \begin{equation*}
        \Big\|\big[\pi^{\varphi}_{p}(b_{ij})\big]_{i,j=1}^{n}\Big\|_{\ell} \leq \Big\|\big[(b_{ij})\big]_{i,j=1}^{n}\Big\|_{\alpha_{p}}, \; \text{for}\; p= 1,\hdots,m.
    \end{equation*}
    It follows from the above arguments that 
    \begin{align}\label{Equation: leqalpha}
      \notag  \Bigg\|(\varphi_{T})_{n}(A_{1}, A_{2}, \hdots, &A_{2m-1}) \begin{bmatrix}h_{1} \\ \vdots\\h_{n} \end{bmatrix} \Bigg\|_{\mathcal{H}_{\ell}^{n}} \\
      & \leq \|A_{1}\|_{\alpha_{m}}\cdots \|A_{m}\|_{\alpha_{1}} \cdots \|A_{2m-1}\|_{\alpha_{m}} \Bigg\| \begin{bmatrix}h_{1} \\ \vdots\\h_{n} \end{bmatrix}\Bigg\|_{\mathcal{H}_{\ell}^{n}}.
    \end{align}
    As $(\Lambda, \leq)$ is a directed set, we can choose $\beta \in \Lambda $ with $ \beta \geq \alpha_{p}$ for each $p = 1, \hdots, m$  and by Equation \eqref{Equation: leqalpha}, we get
    \begin{align}\label{Equation: leqbeta}
        \notag  \Bigg\|(\varphi_{T})_{n}(A_{1}, A_{2}, \hdots, &A_{2m-1}) \begin{bmatrix}h_{1} \\ \vdots\\h_{n} \end{bmatrix} \Bigg\|_{\mathcal{H}_{\ell}^{n}} \\
        & \leq \|A_{1}\|_{\beta}\cdots \|A_{m}\|_{\beta} \cdots \|A_{2m-1}\|_{\beta} \Bigg\| \begin{bmatrix}h_{1} \\ \vdots\\h_{n} \end{bmatrix}\Bigg\|_{\mathcal{H}_{\ell}^{n}}.
    \end{align}
    It is clear from Equation \eqref{Equation: leqbeta} that for every $\ell \in \Omega$, there is a $\beta \in \Lambda$ such that 
    \begin{equation*}
        \big\|(\varphi_{T})_{n}(A_{1}, A_{2}, \hdots, A_{2m-1})\big\|_{{\ell}} \leq 1,\; \text{whenever}\; \|(A_{1}, A_{2}, \hdots, A_{2m-1})\|_{\beta} \leq 1. 
    \end{equation*}
    This implies that $\varphi_{T}$ is local completely contractive. By following the similar arguments, the result holds true for $k = 2m$. We conclude that $\varphi_{T} \in ICC_{loc}^{s}(\mathcal{A}^{k}, C^{\ast}_{\mathcal{E}}(\mathcal{D}))$.
    
 Our next aim is to show that $\varphi_{T}$ is a local completely positive map. Suppose that $k = 2m-1$, $n \in \mathbb{N}$ and $ A_{t} = \big[a_{t_{ij}}\big]_{i,j=1}^{n}\in M_{n}(\mathcal{A})$, for $t = 1, \hdots, 2m-1$. Then for every $\ell \in \Omega$ and $h_{1}, h_{2}, \hdots, h_{n} \in \mathcal{H}_{\ell}$, we have 
    \begin{align*}
        \Bigg\langle (\varphi_{T})_{n}(A_{1}, A_{2}, \hdots, A_{2m-1}) \begin{bmatrix} h_{1}\\ \vdots \\ h_{n} \end{bmatrix}, \; \begin{bmatrix} h_{1}\\ \vdots \\ h_{n} \end{bmatrix} \Bigg\rangle_{\mathcal{H}_{\ell}^{n}} 
        \end{align*}
        \begin{align*}
        &= \sum\limits_{i=1}^{n} \Big\langle\sum\limits_{j=1}^{n} \sum\limits_{r_{1},\hdots,r_{2(m-1)} = 1}^{n}\varphi_{T}(a_{1_{ir_{1}}}, \hdots, a_{{2m-1}_{r_{2(m-1)}j}})h_{j}, \; h_{i} \Big\rangle_{\mathcal{H}_{\ell}}\\
        &= \sum\limits_{i=1}^{n} \Big\langle\sum\limits_{j=1}^{n} \sum\limits_{r_{1},\hdots,r_{2(m-1)} = 1}^{n}V_{\varphi}^{\ast}T\pi^{\varphi}_{1}(a_{m_{r_{m-1}r_{m}}}) \cdots \pi^{\varphi}_{m}(a_{1_{ir_{1}}}a_{{2m-1}_{r_{2(m-1)}j}})V_{\varphi}h_{j}, \; h_{i} \Big\rangle_{\mathcal{H}_{\ell}}\\
        &= \sum\limits_{i=1}^{n} \Big\langle\sum\limits_{j=1}^{n} \sum\limits_{r_{1},\hdots,r_{2(m-1)} = 1}^{n}\pi^{\varphi}_{m}(a_{1_{ir_{1}}}) \cdots \pi^{\varphi}_{1}(a_{m_{r_{m-1}r_{m}}}) \cdots \pi^{\varphi}_{m}(a_{{2m-1}_{r_{2(m-1)}j}})\sqrt{T}V_{\varphi} h_{j}, \; \sqrt{T}V_{\varphi} h_{i} \Big\rangle_{\mathcal{H}^{\varphi}_{\ell}}\\
        &= \Big\langle (\pi^{\varphi}_{1})_{n}(A_{m}) \xi, \; \eta\Big\rangle_{{\mathcal{H}_{\ell}^{\varphi}}^{\oplus n}},
    \end{align*}
    where 
    \begin{equation*}
        \xi = \prod\limits_{p =2}^{m}(\pi^{\varphi}_{p})_{n}(A_{m+p-1}) \begin{bmatrix}
    \sqrt{T}V_{\varphi}h_{1} \\ \vdots \\ \sqrt{T}V_{\varphi}h_{n}
    \end{bmatrix}\;
    \text{and}\;
        \eta = \prod\limits_{p =2}^{m}(\pi^{\varphi}_{p})_{n}(A_{m+1-p}^{\ast})\begin{bmatrix}
    \sqrt{T}V_{\varphi}h_{1} \\ \vdots \\ \sqrt{T}V_{\varphi}h_{n}
    \end{bmatrix}.
    \end{equation*}
    Since $(\pi_{p}^{\varphi})_{n}$ is a local contractive map, there exists an $\alpha_{p} \in \Lambda$ for each $p = 2, \hdots, m$ such that 
    \begin{equation}\label{Equation: foreachalphap}
        \Big\| (\pi_{p}^{\varphi})_{n}(X)\Big\|_{\ell} \leq \|X\|_{\alpha_{p}}.
    \end{equation}
    We can choose $\alpha \geq \alpha_{p}$ for each $p$ and by Equation \eqref{Equation: foreachalphap}, we get 
    \begin{equation}\label{Equation: foreachp}
        \Big\| (\pi_{p}^{\varphi})_{n}(X)\Big\|_{\ell} \leq \|X\|_{\alpha}, \; \text{for every }\; p = 2, \hdots, m.
    \end{equation}
    Also by using the fact that $(\pi^{\varphi}_{1})_{n}$ is a local positive map, there exist a $\beta \in \Lambda$ such that 
    \begin{equation}\label{Equation: pi1}
        (\pi_{1}^{\varphi})_{n}(B) \geq_{\ell} 0, \; \text{whenever}\; B\geq_{\beta} 0.
    \end{equation}
    Since $(\Lambda, \leq )$ is a directed set, we can choose $\gamma \in \Lambda$ with $\alpha, \beta \leq \gamma$.
    Suppose that 
    \begin{align}\label{Eq: gamma}
    (A_{1}, \hdots, A_{m}, \hdots, A_{2m-1}) =_{\gamma}(A_{2m-1}^{\ast}, \hdots,A_{m}^{\ast}, \hdots, A_{1}^{\ast}) \end{align}
    and $A_{m} \geq_{\gamma} 0 $, then due to the relation $\alpha \leq \gamma$, we have
    \begin{equation}\label{Equation: alphasymmetric}
        A_{m+p-1} =_{\alpha} A_{m+1-p}^{\ast} \; \text{and}\; A_{m} \geq_{\alpha}0, \; \text{for each}\; p = 2, \hdots, m. 
    \end{equation}
    It follows from Equations \eqref{Equation: foreachp}, \eqref{Equation: alphasymmetric} that 
    \begin{equation*}
        \Big\| (\pi^{\varphi}_{p})_{n}(A_{m+p-1} - A_{m+1-p}^{\ast}) \Big\|_{\ell} \leq \Big\|A_{m+p-1}-A_{m+1-p}^{\ast} \Big\|_{\alpha} = 0.
    \end{equation*}
    That is, 
    \begin{equation*}
        (\pi_{p}^{\varphi})_{n}(A_{m+p-1})\big|_{{\mathcal{H}_{\ell}^{\varphi}}^{\oplus n}} = (\pi_{p}^{\varphi})_{n}(A_{m+1-p}^{\ast})\big|_{{\mathcal{H}_{\ell}^{\varphi}}^{\oplus n}}, \; \text{for each}\;p = 2, \hdots, m.
    \end{equation*}  
    This implies that $\xi = \eta$ and hence by Equation \eqref{Equation: pi1}, we conclude that
    \begin{equation*}
        \Big\langle (\pi_{1}^{\varphi})_{n}(A_{m}) \xi,\; \xi \Big\rangle_{{\mathcal{H}_{\ell}^{\varphi}}^{\oplus n}} \geq 0 , \; \text{whenever}\; A_{m} \geq_{\beta} 0.
    \end{equation*}
    As $\gamma \geq \beta$, it is clear that $A_{m} \geq_{\gamma} 0$ whenever $A_{m} \geq_{\beta} 0$. Therefore, for each $\ell \in \Omega$, one can choose $\gamma \in \Lambda$ as described in Equation  
    \eqref{Eq: gamma} such that
    \begin{equation*}
        \Bigg\langle (\varphi_{T})_{n}(A_{1}, A_{2}, \hdots, A_{2m-1}) \begin{bmatrix} h_{1}\\ \vdots \\ h_{n} \end{bmatrix}, \; \begin{bmatrix} h_{1}\\ \vdots \\ h_{n} \end{bmatrix} \Bigg\rangle_{\mathcal{H}_{\ell}^{n}} \geq 0,
    \end{equation*}
    whenever $(A_{1}, \hdots, A_{m}, \hdots, A_{2m-1}) =_{\gamma} (A_{2m-1}^{\ast}, \hdots, A_{m}^{\ast},\hdots, A_{1}^{\ast})$ and $A_{m} \geq_{\gamma}0$ and incase $A_{j}=_{\gamma}0$ for some $1\leq j\leq 2m-1,$
    then $(\varphi_{T})_{n}=_{\ell}0$ for all $n.$
    This shows that $\varphi_{T}$ is a local completely positive map.

Suppose that $k = 2m$. Then
\begin{align*}
    \Bigg\langle (\varphi_{T})_{n}(A_{1}, \hdots, A_{2m}) \begin{bmatrix}
    h_{1}\\ \vdots \\ h_{n}
    \end{bmatrix}, \;  \begin{bmatrix}
    h_{1}\\ \vdots \\ h_{n}
    \end{bmatrix} \Bigg\rangle_{\mathcal{H}_{\ell}^{n}} = \Big\langle (\pi_{1}^{\varphi})_{n}(A_{m}A_{m+1})\widetilde{\xi}, \; \widetilde{\eta} \Big\rangle_{{\mathcal{H}_{\ell}^{\varphi}}^{\oplus n}},
\end{align*}
where 
\begin{equation*}
    \widetilde{\xi} = \Big( \prod\limits_{p=2}^{m}(\pi_{p}^{\varphi})_{n}(A_{m+p})\Big) \begin{bmatrix}
    \sqrt{T}V_{\varphi}h_{1}\\ \vdots \\ \sqrt{T}V_{\varphi}h_{n}
    \end{bmatrix}
\; \text{and} \;
    \widetilde{\eta} = \Big( \prod\limits_{p=2}^{m}(\pi_{p}^{\varphi})_{n}(A_{m+1-p}^{\ast})\Big) \begin{bmatrix}
    \sqrt{T}V_{\varphi}h_{1}\\ \vdots \\ \sqrt{T}V_{\varphi}h_{n}
    \end{bmatrix}.
\end{equation*}
Following the similar arguments as in the case of $k=2m-1$ and from Equations \eqref{Equation: foreachp}, \eqref{Equation: pi1}, we conclude that there is a $\gamma \in \Lambda$ with 
\begin{equation*}
    (A_{1}, \hdots, A_{m}, \hdots, A_{2m}) =_{\gamma} (A_{2m}^{\ast}, \hdots, A_{m+1}^{\ast}, \hdots, A_{1}^{\ast}),
\end{equation*}
then $\widetilde{\xi} = \widetilde{\eta}$ and 
\begin{equation*}
    \Big\langle (\pi_{1}^{\varphi})_{n}(A_{m}^{\ast}A_{m})\widetilde{\xi}, \; \widetilde{\xi} \Big\rangle_{{\mathcal{H}_{\ell}^{\varphi}}^{\oplus n}} \geq 0.
\end{equation*}
This implies that $\varphi_{T}$ is a local completely positive. 

\noindent Proof of $(ii):$ We claim that $(\varphi-\varphi_{T}) \in ICPCC_{loc}(\mathcal{A}^{k}, C^{\ast}_{\mathcal{E}}(\mathcal{D}))$. For every $a_{1}, a_{2}, \hdots, a_{k} \in \mathcal{A}$, we have 
    \begin{align*}
        &(\varphi-\varphi_{T})(a_{1},\hdots,a_{m},\hdots, a_{k}) \\
        &=  \left\{ \begin{array}{cc}
             V_{\varphi}^{\ast}(I-T)\pi^{\varphi}_{1}(a_{m}) \pi^{\varphi}_{2}(a_{m-1}a_{m+1})\cdots \pi^{\varphi}_{m}(a_{1}a_{2m-1})V_{\varphi}& \text{if}\; k = 2m-1  \\
             &\\
             V_{\varphi}^{\ast}(I-T)\pi^{\varphi}_{1}(a_{m}a_{m+1}) \pi^{\varphi}_{2}(a_{m-1}a_{m+2})\cdots \pi^{\varphi}_{m}(a_{1}a_{2m})V_{\varphi}& \text{if}\; k = 2m
        \end{array}\right.\\
        &= \varphi_{(I-T)}.
    \end{align*}
    Since $(I-T) \in \bigcap\limits_{p=1}^{m}\pi^{\varphi}_{p}(\mathcal{A})^{\prime}\cap C^{\ast}_{\mathcal{E}^{\varphi}}(\mathcal{D}^{\varphi})$ and $0 \leq (I-T) \leq I$,  it follows from the proof of $(i)$ that  $\varphi_{(I-T)} \in ICPCC_{loc}(\mathcal{A}^{k}, C^{\ast}_{\mathcal{E}}(\mathcal{D})$. Equivalently, we obtain  $\varphi_{T} \leq \varphi$.
\end{proof}
The above results are a suitable recipe for  Radon-Nikod\'ym theorem for local completely positive invariant $k$-linear maps. This will provide a characterization of the set of all maps that are dominated by a fixed map from the class $ICPCC_{loc}(\mathcal{A}^{k}, C^{\ast}_{\mathcal{E}}(\mathcal{D}))$.
\begin{thm}\label{Theorem: Radon-Nikodym}
Let $\varphi, \psi \in ICPCC_{loc}^{s}\big({A}^{k}, C^{\ast}_{\mathcal{E}}(\mathcal{D})\big)$ and $m = \big[\frac{k+1}{2}\big].$  Then $\psi \leq \varphi$  if and only if there is a unique $ \Delta_{\varphi}(\psi) \in \bigcap\limits_{p=1}^{m}\pi_{p}^{\varphi}(\mathcal{A})^{\prime} \cap C^{\ast}_{\mathcal{E}^{\varphi}}(\mathcal{D}^{\varphi})$ with $0 \leq \Delta_{\varphi}(\psi) \leq I$ such that $\psi = \varphi_{\Delta_{\varphi}(\psi)}$.
\end{thm}
\begin{proof}
Since $\psi\in ICPCC_{loc}^{s}(\mathcal{A}^{k}, C^{\ast}_{\mathcal{E}}(\mathcal{D}))$ is a local completely positive map, by Theorem \ref{Theorem: Multilinear Stinespring} we have a minimal Stinespring's triple $\Big( \{\pi_{p}^{\psi}\}_{p=1}^{m}, V_{\psi}, \{\mathcal{H}^{\psi}; \mathcal{E}^{\psi}; \mathcal{D}^{\psi}\}\Big)$ associated to $\psi$.
Suppose that $\psi \leq \varphi$. Then by Lemma \ref{Lemma: Existenceofboundedoperator}, there is a unique contraction $T \in \mathcal{B}(\mathcal{H}^{\varphi}, \mathcal{H}^{\psi})$ such that $TV_{\varphi} = V_{\psi}$ and 
\begin{equation}\label{Equation: commutantofT}
    T\Big(\prod\limits_{p=1}^{m}\pi^{\varphi}_{p}(a_{p})\Big) \subseteq \Big(\prod\limits_{p=1}^{m}\pi^{\psi}_{p}(a_{p})\Big)T, \; \text{for each}\; a_{p} \in \mathcal{A}.
\end{equation} 
Let us define $\Delta_{\varphi}(\psi):= T^{\ast}T$. Then clearly $\Delta_{\varphi}(\psi) \in \mathcal{B}(\mathcal{H}^{\varphi})$ and $0 \leq \Delta_{\varphi}(\psi) \leq I$. Moreover, $\Delta_{\varphi}(\psi)\big|_{\mathcal{H}^{\varphi}_{\ell}} \in \mathcal{B}(\mathcal{H}^{\varphi}_{\ell})$, for each $\ell \in \Omega$. In view of Equation \eqref{Equation: commutantofT}, we compute that \\

\begin{align*}
    \Delta_{\varphi}(\psi)\pi_{p}^{\varphi}(a) = T^{\ast}T\pi_{p}^{\varphi}(a)\prod\limits_{t=1, t\neq p}^{m}\pi^{\varphi}_{t}(1) \subseteq T^{\ast}\pi_{p}^{\psi}(a)\prod\limits_{t=1, t\neq p}^{m}\pi^{\psi}_{t}(1)T
    = \pi_{p}^{\varphi}(a)\Delta_{\varphi}(\psi),
\end{align*}
for every $a \in \mathcal{A}$ and $p \in \{1, \hdots, m\}$. As a result, $\Delta_{\varphi}(\psi) \in \bigcap\limits_{p=1}^{m}\pi_{p}^{\varphi}(\mathcal{A})^{\prime}\cap C^{\ast}_{\mathcal{E}^{\varphi}}(\mathcal{D}^{\varphi})$. Let $a_{1}, \hdots, a_{k} \in \mathcal{A}$ and $h,g \in \mathcal{D}$. Then there is an $\ell \in \Omega$ with $h,g \in \mathcal{H}_{\ell}$ and 
\begin{align*}
  &\Big\langle \varphi_{\Delta_{\varphi}(\psi)}(a_{1}, \hdots, a_{k})h,\; g\Big\rangle_{\mathcal{H}_{\ell}} \\
    =& \left\{ \begin{array}{cc}
        \Big\langle V_{\varphi}^{\ast}\Delta_{\varphi}(\psi)\pi_{1}^{\varphi}(a_{m}) \cdots \pi_{m}^{\varphi}(a_{1}a_{2m-1})V_{\varphi}h,\; g\Big\rangle_{\mathcal{H}_{\ell}} &\; \text{if}\; k = 2m-1  \\
         & \\
         \Big\langle V_{\varphi}^{\ast}\Delta_{\varphi}(\psi)\pi_{1}^{\varphi}(a_{m}a_{m+1}) \cdots \pi_{m}^{\varphi}(a_{1}a_{2m})V_{\varphi}h,\; g\Big\rangle_{\mathcal{H}_{\ell}} &\; \text{if}\; k = 2m 
    \end{array}\right.\\
    =&\left\{ \begin{array}{cc}
        \Big\langle T\pi_{1}^{\varphi}(a_{m}) \cdots \pi_{m}^{\varphi}(a_{1}a_{2m-1})V_{\varphi}h,\; TV_{\varphi} g\Big\rangle_{\mathcal{H}_{\ell}^{\psi}} &\; \text{if}\; k = 2m-1  \\
         & \\
         \Big\langle T\pi_{1}^{\varphi}(a_{m}a_{m+1}) \cdots \pi_{m}^{\varphi}(a_{1}a_{2m})V_{\varphi}h,\; TV_{\varphi}g\Big\rangle_{\mathcal{H}_{\ell}^{\psi}} &\; \text{if}\; k = 2m 
    \end{array}\right.\\
    =&\left\{ \begin{array}{cc}
        \Big\langle \pi_{1}^{\psi}(a_{m}) \cdots \pi_{m}^{\psi}(a_{1}a_{2m-1})V_{\psi}h,\; V_{\psi} g\Big\rangle_{\mathcal{H}_{\ell}^{\psi}} &\; \text{if}\; k = 2m-1  \\
         & \\
         \Big\langle \pi_{1}^{\psi}(a_{m}a_{m+1}) \cdots \pi_{m}^{\psi}(a_{1}a_{2m})V_{\psi}h,\; V_{\psi}g\Big\rangle_{\mathcal{H}_{\ell}^{\psi}} &\; \text{if}\; k = 2m 
    \end{array}\right.\\
    =&\left\{ \begin{array}{cc}
        \Big\langle V_{\psi}^{\ast}\pi_{1}^{\psi}(a_{m}) \cdots \pi_{m}^{\psi}(a_{1}a_{2m-1})V_{\psi}h,\;  g\Big\rangle_{\mathcal{H}_{\ell}^{\psi}} &\; \text{if}\; k = 2m-1  \\
         & \\
         \Big\langle V_{\psi}^{\ast}\pi_{1}^{\psi}(a_{m}a_{m+1}) \cdots \pi_{m}^{\psi}(a_{1}a_{2m})V_{\psi}h,\; g\Big\rangle_{\mathcal{H}_{\ell}^{\psi}} &\; \text{if}\; k = 2m 
    \end{array}\right.\\
    =&\Big\langle \psi(a_{1}, \hdots, a_{k})h,\; g\Big\rangle_{\mathcal{H}_{\ell}}.
\end{align*}
Therefore, $\psi = \varphi_{\Delta_{\varphi}(\psi)}$. Now conversely, suppose that $\psi = \varphi_{\Delta_{\varphi}(\psi)}$, where $\Delta_{\varphi}(\psi) \in \bigcap\limits_{p=1}^{m}\pi_{p}^{\varphi}(\mathcal{A})^{\prime}\cap C^{\ast}_{\mathcal{E}^{\varphi}}(\mathcal{D}^{\varphi})$ with $0 \leq \Delta_{\varphi}(\psi) \leq I$. Then it follows from  Proposition \ref{Proposition: phiT} that  $\psi \in ICPCC^{s}_{loc}(\mathcal{A}^{k}, C^{\ast}_{\mathcal{E}}(\mathcal{D}))$ and $\psi \leq \varphi$. This completes the proof.
\end{proof}
\begin{cor}
Let $\varphi\in ICPCC_{loc}^{s}(\mathcal{A}^{k}, C^{\ast}_{\mathcal{E}}(\mathcal{D}))$ and $m = \big[\frac{k+1}{2}\big]$. Then the map 
\begin{equation*}
    \eta \colon \Big\{T\in \bigcap\limits_{p=1}^{m}\pi_{p}^{\varphi}(\mathcal{A})^{\prime}\cap C^{\ast}_{\mathcal{E}^{\varphi}}(\mathcal{D^{\varphi}}):\; 0 \leq T \leq I\Big\} \longrightarrow \big[0, \varphi\big]
\end{equation*}
defined by $\eta(T) = \varphi_{T}$ is an affine order isomorphism onto $[0, \varphi]$. 
\end{cor}
\begin{proof}
 We know that $\eta$ is injective map by Lemma \ref{Lemma: introductionofphiT} and it is onto by Theorem \ref{Theorem: Radon-Nikodym}. Suppose that $T_{1}, T_{2} \in \bigcap\limits_{p=1}^{m}\pi_{p}^{\varphi}(\mathcal{A})^{\prime}\cap C^{\ast}_{\mathcal{E}^{\varphi}}(\mathcal{D}^{\varphi})$ with $0 \leq T_{1} \leq T_{2} \leq I$ then $0 \leq (T_{2}-T_{1}) \leq I$ and by Proposition \ref{Proposition: phiT}, we have
\begin{equation*}
    \eta(T_{2} - T_{1}) = \varphi_{(T_{2}-T_{1})} \in ICPCC_{loc}^{s}(\mathcal{A}^{k}, C^{\ast}_{\mathcal{E}}(\mathcal{D})).
\end{equation*}
This shows that $\eta(T_{1}) \leq \eta(T_{2})$. Therefore, $\eta$ is an order isomorphism onto $[0, \varphi]$.
\end{proof}

\subsection*{Author Contributions} All authors contributed equally to this work.
\subsection*{Acknowledgements}
The first named author would like to thank his doctoral supervisor Prof. Anil Kumar Karn for introducing the concept of local operator spaces and locally $C^*$-algebras. The second named author thanks NBHM (National Board for Higher Mathematics, India) for financial support with ref No.
[0204/66/2017/R\&D-II/15350] and BITS Pilani. 

Both authors wish to express their sincere gratitude to Prof. B V Rajarama Bhat for mathematical discussions, support, encouragement and are grateful to Indian Statistical Institute, Bangalore for providing the necessary facilities to carry out this work. The authors thank anonymous reviewers for their careful reading and valuable comments.

\end{document}